
\input epsf

\loadeurm

\loadeusm

\magnification=\magstep1

\documentstyle {amsppt}

\def\a{\alpha}
\def\b{\beta}
\def\g{\gamma}
\def\G{\Gamma}
\def\d{\delta}
\def\D{\Delta}
\def\e{\epsilon}

\def\k{\kappa}
\def\l{\lambda}

\def\c{\chi}

\def\p{\pi}

\def\r{\rho}

\def\si{\sigma}
\def\Si{\Sigma}

\def\t{\tau}

\def\o{\omega}
\def\O{\Omega}

\def\tS{\tilde \Si}
\def\8{_{\infty}}

\def\ti{\tilde}
\def\ps{$p_1$-structure\ }
\def\pss{$p_1$-structures\ }

\voffset=-0.5in

\hfuzz 27pt

\topmatter \title Invariants for 1-dimensional cohomology
classes arising from TQFT   \endtitle \rightheadtext{Invariants
arising from TQFT} \author Patrick M. Gilmer \endauthor
 \affil Louisiana State University \endaffil
 \address Department of Mathematics,  Baton Rouge, LA
70803 U.S.A  \endaddress
  \email gilmer\@ math.lsu.edu \endemail

  \abstract Let $(V,Z)$ be a Topological Quantum
Field Theory over a field  $f$  defined on a cobordism category
whose morphisms are oriented  $n+1$-manifolds perhaps
with extra structure (for example a $p_1$  structure and
banded link). Let $(M,\chi)$ be a closed oriented
$n+1$-manifold $M$ with this extra structure together with
$\chi \in  H^1(M).$  Let $M_{\infty}$ denote the infinite cyclic
cover of $M$ given  by $\chi.$    Consider a fundamental
domain $E$ for  the  action of the  integers on  $M_{\infty}$
bounded by lifts of a surface $\Sigma$ dual to   $\chi,$ and in
general position. $E$ can be viewed as a cobordism from
$\Sigma$ to itself.   We give Turaev and Viro's proof of their
theorem  that the similarity class  of the non-nilpotent part
of $Z(E)$ is an invariant.  We give a method to  calculate this
invariant for  the $(V_p,Z_p)$ theories of  Blanchet,
Habegger, Masbaum and  Vogel when $M$ is zero
framed surgery to $S^3$ along a knot K. We give a formula for this invariant
when $K$ is  a
twisted double of  another knot.    We obtain formulas for the
quantum invariants of branched covers of knots, and  unbranched covers of
0-surgery to $S^3$ along knots. We study periodicity among the quantum
invariants of Brieskorn manifolds.  We give an upper bound on the quantum
invariants of branched covers of fibered knots. We also  define finer
invariants for pairs $(M,\chi)$  for TQFT's over
Dedekind domains.  We  use  these ideas to
study isotopy invariants of banded links in $S^1 \times  S^2.$
\endabstract
  \endtopmatter

\document
\centerline{First Version: 1/5/94 \quad This Version: 11/02/95}
\head Introduction\endhead

Witten  conceived of topological quantum field  theory  and
related it to the  Jones polynomial \cite{W}. Axioms, based on
Segal's axioms for  conformal field theory,  were given by
Atiyah \cite{A1,A2}. The first  rigorous development of the
projective version of the TQFT's which most  concern us here
was given by Reshetikhin-Turaev \cite{RT}. There are  now
several rigorous mathematical approaches to topological
quantum field  theory. We have used the work of Blanchet,
Habegger, Masbaum and  Vogel \cite{BHMV1,MV1} as a
foundation for our work, as it is the most complete
development of the  subject for someone with our
background. We have also been influenced by  the papers of
Lickorish \cite{L1}  and   Walker \cite{Wa1} which will  not
be explicitly referred to below. If we have a
$n+1$-manifold fibered over a circle,  and a TQFT in $n+1$ dimensions, then the
monodromy  induces an automorphism of the vector space
associated to the fiber. This construction was generalized to $n+1$-manifolds
$M$ together with a one dimensional cohomology class $\chi$ by Turaev and Viro
\cite{TV} as follows. One considers a fundamental domain $E$ for the
action of the  integers on $M_{\infty},$ the infinite cyclic cover of M. $E$
can be viewed as a cobordism from
a surface to itself.   We give Turaev and Viro's proof of their
theorem  that the similarity class  of the non-nilpotent part
of the induced endomorphism is an invariant.

In \S1, we describe the  Turaev-Viro module of  $(M,\chi)$, which Turaev and
Viro conceived of as being somewhat analogous to the the Alexander module of a
knot, but with a TQFT replacing homology.  We study a number of properties of
the Turaev-Viro module and its associated invariants. In particular, we relate
these invariants to TQFT invariants of the  finite cyclic covers of $M$ given
by $\c.$  In \S2, we show how this invariant may be refined if we are working
with a TQFT defined over a Dedekind domain, rather than a field. The results of
the first two sections are  axiomatic and apply
to any  TQFT and more generally to many linearizations of cobordism categories.
By a linearization over a ring $d$, we mean a functor from a cobordism category
to a category of modules over $d.$ If target category is a category of finitely
generated modules over $d,$ the linearization is said to be finite.
In particular, \S1 and \S2  may be applied to the $Z_p$
theories of  \cite{BHMV1} and the theories of Frohman and
Nicas \cite {FN1,FN2}).
In \S 3, we discuss various issues involving  $V_p$ theories and
$p_1$-structures.     In \S 4 we study banded links $L$ in $S^1 \times S^2$
which are null homologous modulo two. We define a restricted cobordism category
$\eusm C$ and a finite linearization of  $\eusm C$ over $\Bbb Z[A,A^{-1}].$ We
obtain
in this way a polynomial  invariant $D(L) \in \Bbb Z[A,A^{-1}],$ as well as an
invariant of $L$ which is a similarity class of automorphisms of modules over
$\Bbb Q[A,A^{-1},\frac{1}{D(L)}].$ For almost all $p$, these invariants
specialize to invariants associated to $L$ by the $Z_p$ theory as above.
In \S5, we adapt Rolfsen's method \cite{R} of calculating the Alexander module
to the problem of calculating the Turaev-Viro module associated to 0-framed
surgery along a knot. We  study  twisted doubles of knots in detail.
We  also calculate the invariant
for the knot $8_8$ for $p = 5$. We also prove a number of general results.   We
use a result of Casson and Gordon
to obtain a restriction on this invariant for a fibered ribbon knot.

In \S 6, we use these results
to calculate the quantum invariant $< \ >_p$ for the finite cyclic covers of
0-framed surgery along knots. In \S 7, we introduce certain
colored invariants of knots which are necessary to give a good formula for the
Turaev-Viro modules of a connected sum. These same colored invariants are then
used in  \S 8 to give formulas for the quantum invariants of the branched
cyclic covers of knots. We use these  formulas to give closed formulas for $< \
>_5$ for all the branched cyclic covers of the trefoil, the figure eight, the
stevedore's knot, and the untwisted double of the figure eight. In fact, using
our formulas, it is an easy matter to calculate recursively $< \ >_5$ for all
the branched cyclic covers of a twisted double of a knot $J$, once we know two
values of the Kauffman polynomial of $J.$ The same may be said for the
unbranched cyclic covers of zero surgery to $S^3$ along a twisted double of
$J.$ We then derive some periodicity results  for quantum invariants of
Brieskorn manifolds and more generally branched covers of fibered knots with
periodic monodromy.
For fibered knots whose monodromy is not necessarily periodic, we
 obtain an upper bound for their quantum invariants. In \S 9, we consider the
extent to which these invariants are skein invariants. In \S 10, we discuss
$V_p(\Si)$ for $p$ odd and $\Si$ disconnected. In \S 11, we compare when
possible our calculations with other calculations and methods. The afterword
has some final conjectures and other remarks. In the interest of the reader,
we will frequently derive a result, and then subsequently derive a more general
result by a more difficult proof, or a proof requiring more background. We used
Mathematica running on a NeXT computer for our calculations.

 We wish to thank Oleg Viro, Gregor  Masbaum, Pierre Vogel, Neal
Stoltzfus, Rick Litherland, Larry Smolinsky, Steve  Weintraub,
Chuck Livingston, Paul Melvin, Bill Hoffman, Jorge Morales, and Bill
Adkins for useful conversations.

\document

\head \S 1 The Turaev-Viro module   \endhead

Suppose we have a concrete cobordism category $C$ in the
sense of  \cite{BHMV1,(1.A)}. The objects of $C$ are compact
oriented manifolds  of dimension $n$ with perhaps some extra
structure. A morphism from  $\Si$ to $\Si'$ is a   compact
oriented manifold $M$ of dimension $n+1$  with perhaps
some extra structure together with a diffeomorphism of
$\partial M $ to the disjoint union $ -\Si \coprod \Si'$, up to
equivalence. We call such a manifold  a cobordism from $
\Si$  to $\Si'$. Two such cobordisms are equivalent if  there is
a diffeomorphism between them respecting the
diffeomorphism of the  boundary with $ -\Si \coprod \Si'$. In
other words, the relevant diagram  must commute.
Moreover the extra structure on $M$ must induce the  extra
structure on $-\Si \coprod \Si'$. We  assume
codimension-zero  submanifolds in general position inherit
this structure. In addition we  assume the compact
codimension zero submanifolds in general position in  an
infinite  cyclic covering space of a $n+1$-manifold with this
structure  inherit such a structure from their base.  An
especially important example is the category $C^{p_1}_2$
whose objects are closed smooth 2-manifolds  with $p_1$
structure containing a banded (an interval passing through
 each point) collection of
points  and  whose morphisms are smooth 3-manifolds with
$p_1$ structure  containing a  banded link
\cite{BHMV1}.  A banded link in a  3-manifold is  an embedded
oriented surface diffeomorphic to the product  of a
1-manifold with an interval which meets the boundary of the
3-manifold in the product of the boundary of the 1-manifold
with the  interval.

Let $f$ be a field with involution $\l \rightarrow \bar \l$
(perhaps  trivial).  Next suppose we have finite linearization of $C$
over $f.$ It
assigns a $k$-vector space $V(\Si)$  to an object $\Si$ and a
linear transformation  $Z(M)$ to a morphism  $M$. We let
$(V,Z)$ denote this functor.    If we were to follow  the notation  of
\cite{BHMV1,(1.A)},  we should  denote $Z(M)$ by $Z_M$,
and  reserve $Z(M)$ for the morphism induced by the
manifold $M$ viewed as a  cobordism from $\emptyset$ to
$\partial M$.  However this gets more difficult to read as
subscripts proliferate, and it  will be clear from context how
we are thinking of $M$ as a cobordism  from which part of the
boundary to which other part of the boundary.

Let $M\8$ denote the infinite cyclic cover of $M$ classified by
$\chi$, let  $\pi$ denote the projection and $T$ denote the
generating covering  transformation.  Suppose $\g$ is a path
covering a loop on which $\c$  evaluates to 1, then our
convention is that $\g(1) = T(\g(0))$.   Let $\tS$ be any lift of
$\Si $ in $M\8$. Let $E(\tS )$ be the compact  submanifold of
$M\8$ with  boundary $-\tS  \coprod   T\tS $ which may view
as a cobordism from   $\tS$ to $T\tS.$   $E(\tS )$ is a
fundamental domain for the action of the  integers on $M\8$.
If we take $\Si $ to be in general position, $E(\tS )$  defines a
morphism in $C$. Since the projection $\p$ defines a specific
diffeomorphism from any lift of $\Si$ and $\Si$ preserving
any extra  structure, we may regard $Z(E(\tS ))$ as an
endomorphism of $V(\Si)$.   Note that $E(\tS)$ is
diffeomorphic to the exterior of $\Si,$  i.e. M minus  an open
tubular neighborhood of $\Si $.  However this
diffeomorphism  does not necessarily preserve  extra
structure.

Alternatively  we may define $E(\Si)$ to be M ``slit'' along $\Si
.$   This is the $n+1$-manifold obtained by  replacing a
tubular neighborhood of $\Si$ by $\Si \times [-1,0] \coprod \Si
\times [0,1].$  In other words, we take $M$ after we have
replaced each point of $\Si$ by two points and defined
neighborhood systems for these points appropriately. $E(\Si
)$  is a  $n+1$-manifold with structure with boundary $-\Si
\coprod \Si.$

Given a linear endomorphism $\Cal Z$ of a finite dimensional
vector space  $\Cal V$, $\Cal V$ has a canonical $\Cal
Z$-invariant direct sum  decomposition as $\Cal V_0 \oplus
\Cal V_\flat$, where $\Cal Z$ restricted to  $\Cal V_0$ is
nilpotent and $\Cal Z$ restricted to $\Cal V_\flat$ is an
automorphism, denoted $\Cal Z_\flat$. Here $\Cal V_0 =
\cup_{k \ge 1}  \text{Kernel}(\Cal Z^k),$ and $\Cal V_\flat =
\cap_{k \ge 1}  \text{Image}(\Cal Z^k).$   We let $\Cal M (\Cal
Z)$ denote  $\Cal V_\flat$  viewed as a $f[t,t^{-1}]$-module
where $t$ acts by  $\Cal Z_\flat$. We observe that
 one actually has $\Cal V_0 =  \text{Kernel}(\Cal Z^{\dim(\Cal
V)}),$ and $\Cal V_\flat =   \text{Image}(\Cal Z^{\dim(\Cal
V)}).$    Note $\dim(\Cal V_\flat)$ is simply  the number of
nonzero eigenvalues of $\Cal Z$ counted with multiplicity.
$\Cal V_0$ is also known as the generalized 0-eigenspace.

Let $\Cal M_Z(M,\c)$ denote the $k[t,t^{-1}]$-module $\Cal M (Z(E(\Si)))$.
Let $\Cal A_Z(M,\c)$ denote the  automorphism $Z(E(\Si))_\flat.$
In most cases, we will drop the subscript $Z.$

\proclaim{ Theorem (1.1)(Turaev and Viro)} The module $\Cal
M(M,\c)$ is a well defined invariant  of the pair $(M,\c)$ up to
isomorphism. In other words, $\Cal A(M,\c)$ is well defined up
to similarity class.   \endproclaim

\demo{Proof}  Let $\Si$, and $\ti \Si $ be as above, and let  $E=
E(\tS).$  For $k$ an integer, define $\tS_k  = T^{k}\tS.$  For a
negative integer $k<0$, let $E_k = \cup_{k\le i \le -1}  T^{i}E.$ For
a positive integer $k$, let $E_k = \cup_{0\le i \le k-1}  T^{i}E.$     Of
course $\tS_k$  has a natural diffeomorphism with $\Si$. We
can use this diffeomorphism  to give a decomposition
$V(\tS_k) =  V(\tS_k)_0 \oplus V(\tS_k)_\flat.$  If $\Si' $  is a
second oriented surface  dual to $\c$ in general position and,
let $\ti \Si' $ be any lift of $\Si'$ which
 is disjoint from  $\ti \Si $ and  which lies in  $\cup_{i\ge 0}
T^{i}E$. We consider also the case that $\Si'= \Si$, but $\tS'\ne
\tS.$   Let $W$ be the compact  submanifold of $M\8$ with
boundary $-\tS  \coprod  \tS'.$ Define $E'_k,$  $\tS'_k$ and
$V(\tS'_k)_\flat$ analogously to the unprimed items. The
result  will follow from three lemmas:\enddemo

\proclaim {Lemma (1.2)} $Z(W):V( \tS) \rightarrow V( \tS')$
will send  $V( \tS)_\flat$  to $V( \tS')_\flat.$ \endproclaim

\demo{Proof}     Suppose that $x\in V( \tS)_\flat,$ and $Z(W)x
= x'.$ For every $k<0$, there is a $y \in V(\tS_k)$ such that
$Z(E_k)y=x.$ Let $Z(T^{- k}W)y = y'.$ Since $T^{-k}W \cup E'_k
=  E_k \cup W$, we have that  $Z(E'_k)y'=y.$ Thus we have
proved that $Z(W):V( \tS) \rightarrow V( \tS')$ will send  $V(
\tS)_\flat$ to $V( \tS')_\flat.$ We will denote this homomorphism
$Z(W)_\flat.$
\qed\enddemo

\proclaim {Lemma (1.3)} $Z(W)_\flat$ sends $V( \tS)_\flat$
isomorphically to  $V( \tS')_\flat$ \endproclaim

\demo{Proof} For some negative $k$, $\tS'_k$ will lie in
$\cup_{i<0}T^i  (E(\Si ))$.   Let $X$ be the compact submanifold
of $M\8$ with boundary  $-\tS'_k  \coprod  \tS,$ then $X\cup
W = E'_k.$  By functoriality, $Z(W)_\flat  \circ Z(X)_\flat = Z
(E'_k)_\flat.$ $Z(E'_k)$ may be   naturally identified with
$Z(E')^k$, which is an isomorphism when  restricted to $
V(\Si')_\flat$.  Thus $Z(E'_k)_\flat$ is an isomorphism.  Thus
$Z(W)_\flat$ is surjective.

Similarly $W\cup T^{-k}X = E_{-k},$ and $Z(E_{-k})_\flat$ is an
isomorphism. By functoriality, $Z(T^{-k}X)_\flat  \circ
Z(W)_\flat  = Z (E_k)_\flat,$  thus $Z(W)_\flat$ is injective.
 \qed\enddemo

\proclaim{ Lemma (1.4)}Assume $\ti \Si' $ is a lift of $\Si' $
which lies  in $\cup_{i\ge 1}T^i (E(\Si )).$    Let
 $U$ be the compact  submanifold of $M\8$ with boundary $-T
\tS   \coprod  \tS' .$  $Z(U)_\flat\circ Z(\tilde E)_\flat=
Z(E'_{-1})_\flat \circ Z(T^{-1}U)_\flat.$ Thus $Z(E)_\flat $ is
similar to $Z(E')_\flat .$  \endproclaim

\demo{Proof}  $ \tilde E \cup U= t^{-1}U \cup E'_{-1} .$
\qed\enddemo

\noindent {\bf Remarks.} Walker earlier observed that
 $\text{rank} ( Z(\tilde E)_\flat)$\  ($=\dim_f(\Cal M  (M,\c) )$)
is an invariant \cite{Wa2}. Suppose $\dim V(\Si)$ for $\Si$
connected depends only on the genus of $\Si$.  Suppose that
$\dim V(\Si)$ is a increasing  function of the genus $\Si$, as is
true for  Witten's TQFT's. Then $\dim_f(\Cal M(M,\c))$ can be
used to give  lower bounds on the least genus of an embedded
surface dual  to $\c. $ This is just the Thurston norm on
$H_2(M)$ \cite{T}.  Walker discussed  this application.  Turaev
and Viro then strengthened Walker's work.

Let  $\G_Z(M,\c)$ denote the characteristic polynomial of $\Cal
A(M,\c)$. We define  the normalized characteristic polynomial of
a matrix or endomorphism of a free module to be
the characteristic polynomial in $x$ divided by the highest
power of
$x$ dividing this polynomial.   $\G_Z(M,\c)$ is then the normalized
characteristic polynomial of $Z(E)$ for any choice of $\Si$ dual to $\chi$ in
$M.$
It will be convenient to let $D(M,\c)$ denote the constant term of
$\G_Z(M,\c).$

Our convention is that
 the characteristic polynomial of an endomorphism of a zero
dimensional vector space is the constant $1$. We have
$\text{deg}(\G_Z(M,\c))\allowmathbreak = \dim_f(\Cal M((M,\c))).$
$\G_Z(M,\c)$
is analogous to the Alexander polynomial, and should be
called the Turaev-Viro polynomial. A complete set of
invariants  for $\Cal M(M,\c)$ is of course given by the
invariant factors of  $\Cal A(M,\c).$ These are in turn
determined by certain determinantal divisors
\cite{AW,p.312} analogous to the higher Alexander
polynomials.  Two matrices are similar over $f$ if and only if
they are  similar over a larger field
\cite{AW,p.315}.  Thus no information is lost if we
extend our scalars to a larger field. Thus another complete
invariant  would be the Jordan form of $\Cal A(M,\c)$ over
the algebraic closure of $f.$

{\it For the rest of this section,  we suppose that $(V,Z)$ is a
cobordism generated quantization
 i.e. $(V,Z)$ satisfies axioms Q1  Q2 and CG of \cite{BHMV1}}. Given a
vector space $\Cal V$ over $f$, then  $\Cal V^*$ denotes the
vector space with the same underlying Abelian group but
with $\l  v\in \Cal V^*$ given by $\bar \l v \in \Cal V$ for $\l
\in f$ and $v \in \Cal V.$  If $\Cal M$ is a module over
$f[t,t^{-1}]$, then $\Cal M^*$ denotes the conjugate of the
underlying vector space, with $t$ acting the same as before
on the underlying Abelian group. It is clear that a matrix for
the action of $t$  on $\Cal  M^*$, is given by the conjugate of a
matrix for the action of $t$ on $\Cal M.$ Making use of the fact
that a matrix is similar to its transpose, we have the following
proposition:

\proclaim {Proposition (1.5)}  We have:
 $$\Cal M (M,-\c) = \Cal M (M,\c)$$
 $$\Cal M (-M,\c) = \Cal M(M,\c)^*.$$ \endproclaim

\proclaim {Proposition (1.6)} Suppose $M$ is a fiber bundle
over circle with fiber $\Si$ and  let $\c \in  H^1(M)$ be the
cohomology class which is classified by the projection,  then \hfill
$D(M,\c) \overline{D(M,\c)}=1$.  If $z$ is an
eigenvector of  $\Cal A(M,\c)$ with eigenvalue $\l$, either
$<z,z>_{\Si} = 0,$ or $\l\bar \l= 1.$  Moreover eigenvectors with
distinct eigenvalues are orthogonal with  respect $<\ ,\
>_{\Si}.$ In particular if the inner product is definite, then all
the roots of $\G(M,\c)$ have their conjugates as
reciprocals.\endproclaim

\demo{Proof} The monodromy $T$  of this bundle is a
diffeomorphism of  $\Si$ which preserves the \ps  which $\Si$
inherits from $M$. $E(\Si)\in  \tilde \Cal M$ is the mapping
cylinder of $T$. One may check that $T$ is  an isometry of $<\ ,\
>_{\Si}$. Thus the norm of the determinant of a  matrix which
represents this map is one. Then $<z,z>_{\Si}  =<Tz,Tz>_{\Si} = \l
\bar \l <z,z>_{\Si}.$ Similarly if  $z_1$ and $z_2$ are
eigenvectors with distinct eigenvalues $\l_1$ and $\l_2$,
$<z_1,z_2>_{\Si} =<Tz_1,Tz_2>_{\Si} = \l_1 \bar \l_2
<z_1,z_2>_{\Si}.$ \qed \enddemo

{\it For the rest of this section, we suppose that $(Z,V)$
satisfies all the axioms for a TQFT in the sense of
\cite{BHMV1,(1.A)}.} Let $\c_d$ denote $\c$ modulo $d$ and $(M,\c)_d$ denote
the
d-fold cyclic
cover of $M$  classified by $\c_d$ with the induced
structure. By the trace formula of TQFTs   \cite{BHVM, (1.2)},
we have:

\proclaim {Proposition (1.7)} We have:   $Z((M,\c)_d) =
\text{Trace} (\Cal A((M,\c))^d).$\endproclaim

\proclaim {Corollary (1.8)}  $Z((M,\c)_d)$ may be computed
recursively with recursion relation given by $\G(M,\c).$  If $f$
has characteristic zero, then the values of $Z((M,\c)_d)$ for all
$d$, determine $\G(M,\c).$ \endproclaim

\demo{Proof} The first statement just follows from the
Cayley-Hamilton Theorem applied to $\Cal A(M,\c),$ or it may be
obtained from Newton's formula as below.  The trace of  $\Cal
A(M,\c))^d$ is the sum of the dth powers of the eigenvalues of $\Cal
A(M,\c)$ counted with multiplicity. Moreover the coefficients of
$\G(M,\c)$ are, up to sign, the elementary symmetric functions of
the eigenvalues of $\Cal A(M,\c)$ counted with multiplicity. Thus
the initial terms of the sequence $Z((M,\c)_d)$  may also be
computed from the coefficients of $\G(M,\c)$ using Newton's
formula \cite{Ms,  (Problem 16-A)}.   Over a field of characteristic
zero, Newton's formula allows us to calculate the coefficients of
$\G(M,\c)$ recursively from the values of $Z((M,\c)_d).$
\qed \enddemo

\noindent {\bf Remarks.} For example: if $\G (M,\c) = x^2 -
\si_1 x + \si_2$ then $Z(M) = \si_1,$ $Z((M,\c)_2) = \si_1^2 -
2\si_2,$ and $Z((M,\c)_d) = \si_1 Z((M,\c)_{d-1}) -\si_2
Z((M,\c)_{d-2}$ for  $d >2.$ Girard's Formula \cite{MS,
(Problem 16-A)} gives a closed formula for $Z((M,\c)_d)$ in
terms of the coefficients of $\G(M,\c).$ Thus $\G(M,\c)$
contains the same information as the sequence $Z((M,\c)_d).$
However $\G(M,\c)$ is a compact way of organizing this
information.
 If $\G(M,\c)$ has distinct roots, then it determines the
similarity class of  $\Cal A(M,\c)$ and thus the isomorphism
class of $(\Cal M,\c).$

\proclaim {Proposition (1.9)} Suppose $(V,Z)$ is a TQFT defined
on $C_2^{p_1}$ and it satisfies the surgery axiom (S1) of
   \cite{BHMV1}, then we have: $$\Cal A(M \# M',\c\oplus \c') =  \eta
(\Cal A(M,  \c) \otimes \Cal A(M',   \c')).$$ \endproclaim

 \head \S 2 TQFT's over Dedekind Domains \endhead

For each positive integer $p$, \cite{BHMV1} defined cobordism
generated quantizations $(V_p,Z_p)$ which take values in free
finitely generated $k_p$ modules and homomorphism of $k_p$
modules, where $$k_p=\Bbb Z
[\frac{1}{d},A,\k]/(\varphi_{2p}(A),\k^6-u)$$ where
$\varphi_{2p}(A)$ is the $2p$-cyclotomic polynomial in the
indeterminate $A$,  $$d= \cases p,&\text{for $p \neq $3,4,6}\\
	1,&\text{for $p =$ 3,4}\\
 	2,&\text{for $p =6$,}\endcases \text{ and }  u= \cases
A^{-6-\frac{p(p+1)}{2}},&\text{for $p \neq $ 1,2}\\
	1,&\text{for $p =1$}\\
 	A,&\text{for $p =2$.}\endcases $$
We let $A_p$ denote the image of $A$ in $k_p.$ Note that $u_p$,  the image of
$u$ in $k_p,$ is also $1$ for $p= 3$ or $4.$

In every case $k_p$ is the ring of integers of a cyclotomic number
field localized with respect to the multiplicative subset $\{ d^n \ | n
\in \Bbb Z, \ n \ge 0\}$ and so is a Dedekind domain by
\cite{La,p.21} for instance.
 For this reason, we consider now  invariants which may be
defined in this,  and in even more general circumstances which are
potentially stronger than those obtained by passing to the field of
fractions and applying the previous section. {\it In this section we
will assume that $(V,Z)$ is a finite linearization over a Dedekind domain $k.$}

Let $\Cal V$ be a finitely generated module over a Dedekind domain
$k$, $\Cal Z$ a endomorphism of $\Cal V$. Let  $\Cal V_0 = \cup_{k
\ge 1}  \text{Kernel}(\Cal Z^k).$   Let $\Cal V_\sharp=\Cal V /\Cal
V_0.$   Since $\Cal V_0$ is $\Cal Z$ invariant, there is an induced
endomorphism $\Cal Z_\sharp$ of  $\Cal V_\sharp.$ It is clear that
$\Cal Z_\sharp$ is injective.   Let $coker(\Cal Z_\sharp)$ denote the
cokernel of $\Cal Z_\sharp.$
 Because $\Cal Z_\sharp$ is injective, $coker(\Cal Z_\sharp)$ is a
torsion module. To see this consider  the map induced
on  $\Cal V$ modulo its torsion submodule. This map is an
isomorphism after tensoring with the field of fractions of $k$. It
follows that
$\text{coker}(\Cal Z_\sharp)$ is a torsion module. It must also be
finitely generated as $\Cal V$ maps onto it.  Thus  $\text{coker}(\Cal
Z_\sharp)$ is a direct sum of cyclic modules
of the form $k /\text{ann}$ where ann is a non-trivial ideal of $k$
\cite{J,Thm. 10.15}.

 A k-module is of finite length if and only if it is finitely generated.
This follows from the classification of finitely generated  modules
over a Dedekind domain given in \cite{J, Chapt.10}.  Given
$k$-module  $F$ of finite length,  Serre \cite{S,p.14} defined an ideal
of $k$, denoted $\chi_k(F).$  Given a short exact sequence, the value
of $\chi_k$ on the middle term is the product of its value on the
side terms.  $\chi_k(G),$ for $G$  a torsion module, is a nontrivial
ideal.  In fact if
$$G =\bigoplus_{\goth p ,i} \goth p^{e_{i,
\goth p}}$$
then
$$\chi_k(G) =\prod_{ \goth p}\goth p^{(\sum_i
e_{i,\goth p})},$$
where $\goth p$ ranges over all prime ideals of $k,$
and almost all $e_{i,\goth p}$ are zero.  If $k$ is a PID,
then $\chi_k(F)$ is  the order ideal of $F$ as defined by Milnor
\cite{M}. If $G$ is the cokernel of a 1-1 map
of free $k$ modules of finite rank, $\chi_k(G)$ has a nice interpretation.
We need the following result which is less general than \cite{Ba,p. 500}. We
include a proof for the convenience of the reader.

\proclaim{Proposition (2.1)} If $G$ is the cokernel of
$\Cal Z: k^n \rightarrow k^n$ and $\det(\Cal Z)\neq 0$, then
$\chi_k(G)$ is the principal ideal generated by $\det(\Cal Z).$
 \endproclaim

\demo{Proof} We need to see that for each $\goth p$, $\nu_{\goth p}(\det u)=
\sum_i
e_{i,\goth p}.$
 Clearly $\det (u)_\goth p = \det (u_\goth p)$. Also
$\text {cokernel}( u_{\goth p})= \oplus_{i} \goth p^{e_{i,
\goth p}}.$ Thus
$\chi_{k_{\goth p}}
(\text {cokernel}
( u_{\goth p})) = \goth p^{\sum_i e_{i,\goth p}}.$
Finally as $k_{\goth p}$ is a PID and the result to be proved
is  true  for PID's \cite{S,p.17 Lemma 3}, $\det (u_\goth p)=
\chi_{k_{\goth p}}
(\text {cokernel}
( u_{\goth p})). $
 \qed \enddemo

Let $\Cal I (\Cal Z)$ denote $\chi_k(\text{coker}(\Cal Z_\sharp)).$
Let $f$ denote the field of fraction of $k$, and let $k_{\Cal I}$ denote
ring $\{ x \in f\ |\  \nu_\goth p (x) \ge 0\  \forall \goth p \nmid \Cal
I \}.$  Let $\Cal Z_\natural$ denote the endomorphism $\Cal
Z_\sharp \otimes \text{id}_{k_{\Cal I(\Cal Z)}}$,  of  $\Cal V_\sharp
\otimes k_{\Cal I(\Cal Z)},$  which we denote  $\Cal V_\natural.$  Then $\Cal
Z_\natural$ is an isomorphism. Here we have localized as little as
possible such that $\Cal Z_\natural$ is an isomorphism. Let
$\goth M(\Cal Z)$ denote the $k_{\Cal I(\Cal Z)}[t,t^{-1}]$-module given by
the action of  $\Cal Z_\natural$ on $\Cal V_\natural.$

An $n \times n$ matrix $H$ with coefficients in $k$ defines an endomorphism
$\Cal Z_{H}$ of $k^n.$ We let  $H_\natural$ denote
the similarity class of the induced automorphism
$\left( \Cal Z_{H} \right)_\natural$ of a $k_{\Cal I(\Cal Z_{H})}$-module.

\proclaim{Proposition (2.2)}  Let $\Cal  Z$ and $\Cal  Z'$ be two
endomorphisms defined on the finitely generated $k$-modules
$\Cal  V$ and  $\Cal  V'$. Suppose that the induced injections  $\Cal
Z_\sharp$, and $\Cal  Z'_\sharp$ fit into commutative diagram with
$\a$ injective:

$$\CD
 \Cal V_\sharp @>\Cal Z_\sharp>> \Cal  V_\sharp \\
 @V{\a}VV         @VV{\a}V  \\
 \Cal V'_\sharp @>\Cal  Z'_\sharp>> \Cal  V'_\sharp \endCD $$

Then $ \chi_k(\text{coker}(Z_\sharp) ) =
\chi_k(\text{coker}(Z'_\sharp) ).$ \endproclaim

\demo{Proof} We form a short exact sequence of chain complexes all
of whose nonzero terms are concentrated in two dimensions:

$$\CD
 0@>>> \Cal  V_\sharp @>\Cal  Z_\sharp>>  \Cal  V_\sharp@>>>
\text{coker}(\Cal  Z_\sharp) @>>>0\\
 @. @V{\a}VV         @V{\a}VV  @V{\b}VV @.\\
 0@>>>\Cal  V'_\sharp @>\Cal  Z_\sharp>>\Cal  V'_\sharp @>>>
\text{coker}(\Cal  Z'_\sharp)@>>>0\endCD $$

Then the induced long exact sequence of homology is:

$$0 \rightarrow \text{ker}(\b) \rightarrow \text{coker}(\a)
\rightarrow \text{coker}(\a)\rightarrow \text{coker}(\b)
\rightarrow 0.$$

On the other hand we have:

$$0 \rightarrow \text{ker}(\b) \rightarrow \text{coker}(\Cal
Z_\sharp) \rightarrow \text{coker}(\Cal  Z'_\sharp) \rightarrow
\text{coker}(\b) \rightarrow 0.$$

Because of the multiplicative property of $\chi_k$, the result
follows.\qed\enddemo

Let $\goth I_Z(M,\c)$ denote the ideal  $\Cal I(Z(E).$ Let $\goth
M_Z(M,\c)$ denote the module $\goth M (Z(E(\Si)))$. Let $\goth
A_Z(M,\c)$ denote the  automorphism $Z(E(\Si))_\natural.$  In
most cases, we will drop the subscript $Z.$

\proclaim{Theorem (2.3)} $\goth I(M,\c)$, and  the isomorphism
class of $\goth M(M,\c)$ (or the similarity class of $\goth A(M,\c)$ )
are
 invariants  of $(M,\c).$ \endproclaim

\demo{Proof} We use the notations at the beginning of the proof of
Theorem (1.1). Let $V(\Si_k)_\sharp
=V(\Si_k)/V(\Si_k)_0.$ Analogous to Lemma(1.2) we
have:\enddemo

\proclaim {Lemma(2.4)} $Z(W):V( \tS) \rightarrow V( \tS')$ will
send  $V( \tS)_0$  to $V( \tS')_0.$ Thus there is an induced map
 $V( \tS)_\sharp \rightarrow V( \tS')_\sharp $ which we will
denote $Z(W)_\sharp.$ \endproclaim

\demo{Proof}     Suppose that $x\in V( \tS)_0,$ then for some $k>0$,
$Z(E_k)x=0.$  Since $W \cup E'_k= E_k \cup T^k W,$
$$Z(E_k')\circ Z(W) = Z(T^k W)\circ Z(E_k).\tag 2.5$$
 It follows that
 $Z(E')^k(Z(W)x)=0.$ This means $Z(W)x \in V(\tS')_0.$ \qed \enddemo

\proclaim {Lemma(2.6)}
 $Z(W)_\sharp:V( \tS)_\sharp \rightarrow V( \tS')_\sharp $ is
injective. \endproclaim

\demo{Proof} Let  $x \in V( \tS),$ and let $[x]$ denote the
the image of $x$ in $V( \tS)_\sharp.$  Suppose $Z(W)_\sharp ([x])
=0.$  Then $Z(W) (x) \in V(\tS')_0.$ So for some $k$,
$Z(E_k')((Z(W)x)=0.$ For some $l>0,$ $\tS_l \notin W \cup E_k'$.  Let
 $X$  be the compact  submanifold of $M\8$ with boundary $-
\tS'_k   \coprod  \tS_l.$  Thus $Z(X) \circ Z(E_k') \circ Z(W)
(x)=0.$ Since $E_l = W \cup E_k'\cup X$, $x \in V( \tS)_0,$ thus
$[x]=0.$ \qed \enddemo

By Lemma (2.6) and Equation (2.5), the hypothesis of Proposition
(2.2) are now satisifed with $\Cal Z =Z(E),$ and $\Cal Z' =Z(E').$
Thus $\Cal I(Z(E)) =\Cal I(Z(E')).$  Thus $Z(E)_\natural$ is an
isomorphism.  Now as in the proof of
Lemma (1.4), we have that the similarity class of $Z(E)_\natural$
does not depend on the choice of $\Si.$ \qed

Given a finite linearization $(V,Z)$ over a Dedekind domain $k$ as above,
we may of course obtain a finite linearization, say, $(\hat V, \hat Z)$ over
$\hat k,$ the field of fractions of $k,$ by tensoring with $\hat k$. Then we
have
$\goth A_Z(M,\c) \otimes \hat k = \Cal
A_{\hat Z}(M,\c)$
and
$\goth
M_Z(M,\c) \otimes \hat k = \Cal
M_{\hat Z}(M,\c).$
We also define $\G_Z(M,\c)=\G_{\hat Z}(M,\c)$ and
$D_Z(M,\c)=D_{\hat Z}(M,\c).$ Using (2.1), we have:

\proclaim{Proposition (2.7)} Suppose  there is a surface $\Si$ dual to $\c$
such that  $V(\Si)$ is a free $k$ module and $Z(E(\Si))$ is injective, then
$\Cal I (M,\c)$ is the principle ideal generated by $D_Z(M,\c).$ \endproclaim

The hypothesis of (2.7) usually holds in the examples we have studied.
However it does not always hold.

\head \S 3 The  $(V_p,Z_p)$ Theories \endhead

{}From now on, we will mainly be discussing
the  $(V_p,Z_p)$ theories of \cite{BHMV1}. These are defined on
 the cobordism category $C_2^{p_1}$ and on
the larger category $C_{2,q}^{p_1,c}$ \cite{BHMV1,4.6} whose
objects are surfaces with \ps with $q$-colored banded points, and
whose morphisms are 3-manifolds with \ps with a banded
trivalent $q$-colored graph with admissible $q$-coloring. Here a
$q$-coloring assigns to each  edge or framed point an integer from
zero to $q-1$. Here $q$ is $\frac{p-2}{2}$ if $p\ge 4$ and is even,  and is
$p-1$ if $p\ge 3$ and is odd. In the
case $p$ is one of two, we take $q$ to equal two and assume that the union of
the edges of the graph weighted one is a link.
We will call one of these integers a $q$-color. A good $q$-color is simply
a $q$-color if $p$ is even and is an even $q$-color if $p$ is odd.
Recall \cite{BHMV1} a triple of $q$-colors $(i,j,k)$ is called admissible if
$i+j+k \equiv 0 \pmod{2},$
and $i \le j+k,$ $j \le i+k,$ and $k \le i+k.$ We will say that an admissible
triple $(i,j,k)$ is small if in addition $i+j+k < 2q.$ A coloring is said to be
admissible if the colors of the edges meeting at any vertex of order three
form an admissible triple. A coloring is small if these admissible triples are
small.
{}From now on we
will refer to an  admissibly $q$-colored trivalent banded graph as
simply a colored graph.  Note that the notion of a colored graph
includes the notion of a colored link, and plain banded link as
special cases.  A colored graph which happens to be a link (i.e. there
are no 3-valent vertices)  will be called a colored link.

 One may
regard  $C_2^{p_1}$ as a subcategory of $C_{2,q}^{p_1,c}$ by
assigning one uniformly.  The target category for $(V_p,Z_p)$ is the category
of
free finitely generated $k_p$ modules, and module
homomorphisms. The functor from $C_2^{p_1}$ is the
composition of inclusion and the functor from $C_{2,q}^{p_1,c}.$ Thus
we do not really need to distinguish between these two
linearizations in our notation. If the structure of  $M$ contains a
banded link (colored graph) we will denote it by $L$ ($G$).

We say $L$ is an even link in $(M,\c)$ if  $\c$ reduced modulo two
is trivial on the nonoriented fundamental class of $L$.
Otherwise $L$ is an odd link  in
$(M,\c)$. A colored graph is odd or even
according to whether its expansion \cite {BHMV1} is.

\proclaim{Proposition 3.1} If the colored graph in
$(M,\c)$ is odd, $\goth M_{Z_p}(M,\c)=0.$ \endproclaim

\demo{Proof} In this case $V_p(\Si)$ is zero as  a surface with an
odd number of points is not a boundary.\qed \enddemo

\proclaim{Proposition (3.2)}  $\goth A
_{Z_p}(M,\c)$ does not  change if we vary the \ps on $M$ by a
homotopy.  \endproclaim

\demo{Proof}  Suppose we change the $p_1$-structure on $M$
by a homotopy, then we   change the $p_1$-structure on $E$
by a homotopy during which the  $p_1$-structure induced
on the two copies of $\Si$ is identical.  Let
$E'$ denote $E$ with the new $p_1$-structure, and let $\Si'$  denote
$\Si$ with the new $p_1$-structure.   We  use this homotopy of \ps
restricted to $\Si$ to put  a $p_1$-structure on $I \times \Si$.
 As $L\cap \Si$
or $G \cap \Si$ defines some framed  points in $\Si$, $ (-1,1)$ times
these framed points defines a banded link in  $I \times \Si$. Let $P$
denote $I \times  \Si $ equipped  with the
above $p_1$-structure, and banded link.  Let $E''=P \cup E
\cup - P$ glued along the two copies of $\Si$.  $E'$ and
$E''$ both  define morphisms from $\Si'$ to $\Si'$.  There is a
diffeomorphism from  $E'$ to $E'',$  and if we pull back the
$p_1$-structure on $E''$ to $E'$, it  is homotopic to the
$p_1$-structure on $E'$. Similarly if we pull back the  banded
link in $E''$ to $E'$, it is isotopic to the link in $E'$.  Thus our
morphisms $Z_p(E'')$ to $Z_p(E')$ are equal. Clearly $Z_p(E'')$
is  similar to $Z_p(E)$.\qed \enddemo

\proclaim{Proposition (3.3)}   $\k^{-\si (\a(M))}\goth A _{Z_p}(M,\c)$
is invariant as we vary  the \ps on $M.$ \endproclaim

\demo{Proof} If we change the homotopy class of the
$p_1$-structure we  may do that in a small ball neighborhood
well away from $\Si$. This lifts  to a change in the
$p_1$-structure of $E$ which takes place in a small ball  in the
interior of $E$.  Using the functorial properties of $(Z_p,V_p)$,
this will change $Z_p(E)$  and $Z_p(M)$ by the same nonzero
factor.   By \cite{BHMV1,(1.8)}, this  factor is compensated
for by $\k^{-\si(\a)}.$\qed \enddemo

In view of the above, we let $Z_p(M,\c)$ denote $\k^{-\si
(\a)}\goth A  _{Z_p}(M,\c)$.  We also let  $\hat Z_p(M,\c)$ denote $\k^{-\si
(\a)}\Cal A  _{\hat Z_p}(M,\c)$. In this  way, we may remove the
dependence of our invariants on the $p_1$-structure.  The  dependence on
the banding  of the link $L$ or colored graph $G$ is similar. However there is
no  integer invariant of this banding that can be defined in this
generality to  play the role of $\si.$
Since $e_i\in V_p(S^1 \times S^1)$ is an eigenvector for the twist map with
eigenvalue $\mu(s) = (-1)^s A^{s^2+2s}$ \cite{BHVM1,(5.8)}, we can show

\proclaim{Proposition (3.4)}   $Z_p(M,\c)$  is multiplied by
$\mu(s)$  when we change the banding on the colored graph $G$  by adding a
single positive full twist to a single edge colored $s$.
\endproclaim

 The behavior of the $\si$-invariant of \pss  under a cover is related to
signature defects \cite{H} \cite{KM2}:

\proclaim{Proposition (3.5)} Given  $(M,\c)$ with \ps $\a(M)$ and $\c \in
H^1(M)$, we have $\si(\a(M_d)) = d \si(\a(M))- 3\ \text{def\ }(M,\c_d).$
\endproclaim

If $N$ is a morphism from $\emptyset$ to $\emptyset$ in
$C_2^{p_1}$ or $C_{2,q}^{p_1,c},$ we follow \cite{BHMV1} and denote $Z_p(N)\in
k_p$ by $<N>_p.$

Let $\eurm i$ denote an embedding of $k_p$ in $\Bbb C$ which sends $A$ to
$e^{\frac{ \pi i}{p}}$ and sends $\eta$ to a positive  number. There is such an
embedding since $\eurm i(\k^3)$ is only determined up to sign by the choice of
$\eurm i(A).$ Then $[n]= \frac{A^{2n} -A^{-2n}}{A^{2} -A^{-2}}$ will be sent to
 a positive number for $ n< p.$ If $(i_1,i_2,i_3)$ is a small admissible triple
of q-colors, with associated internal colors $(\a,\b, \g)$ then $\eurm
i([k])>0$ if $k$ is one of $i_1,$ $i_2,$ $i_3,$ $\a,$ $\b,$ $\g,$ $\a+\b+\g+1.$
Thus, in this situation, $ (-1)^{\a+\b+\g} \eurm i(<i_1,i_2,i_3>) >0.$
Also for any q-color $c$, $\eurm i([c]) >0.$ \cite{BHMV1,(4.11),(4.14)}
describes a basis for $V_p(\Si)$ given by a small admissible  coloring of a
trivalent graph in a handlebody with boundary $\Si.$ They also describe
the Hermitian form $<\ >_{\Si}$ on $V_p(\Si).$

\proclaim{Proposition (3.6)} If we extend our coefficients to $\Bbb C$ by
$\eurm
i,$ then the form $<\ >_{\Si}$ on $V_p(\Si)\otimes \Bbb C$ is positive
definite.
\endproclaim
 By (1.6) we have:
\proclaim {Corollary (3.7)} Suppose $M$ is a fiber bundle
over a circle with fiber $\Si$ and  let $\c \in  H^1(M)$ be the
cohomology class which is classified by the projection.  Then
the roots of $\eurm i(\G_p(M,\c))$  lie on the unit circle.\endproclaim

By the triangle inequality and (1.7), we have:
\proclaim {Corollary (3.8)} Suppose $M$ is a fiber bundle
over a circle with fiber $\Si$  and  let $\c \in  H^1(M)$ be the
cohomology class which is classified by the projection.  For all
$d,$
$|\eurm i (<(M,\c)_d>_p)| \le \dim V_p(\Si).$ \endproclaim

\proclaim{Proposition (3.9)} Let $\Si$ denote a surface and suppose $\Si \times
S^1$ is given a \ps with $\si$ zero. For $p\ge 3,$ $<\Si \times S^1>_p =
\text{rank}_{k_p} V_p(\Si).$ \endproclaim

\demo{Proof}
 Give $\Si$ a $p_1$-structure, give $S^1$ the \ps coming from  a framing on
$S^1$.
The mapping torus of the identity map on $\Si$, is $\Si \times S^1$ with
the product \ps which we denote by $\a$. $\Si \times S^1$ is the boundary of
$\Si \times D^2,$  and the \ps on $\Si \times S^1$ extends over this 4-manifold
as any \ps on $S^1$ extends to one on $D^2.$ Since the signature of $\Si \times
D^2$ is zero,
$\si(\a)$ is zero. Also $\text{rank}_{k_p} V_p(\Si)$ is the trace of the
identity on
$V_p(\Si).$ \qed \enddemo

 Of course the above proposition is well known except possibly for  nailing
down the $\si$ invariant.

\proclaim {Corollary (3.10)} Suppose $M$ is a fiber bundle
over a circle with fiber $\Si$  with mondromy of period $s.$ Suppose the
colored graph in $M$ is empty. Let $\c \in  H^1(M)$ be the
cohomology class which is classified by the projection. Assume $p\ge 3.$
If $p \equiv 0\pmod{4}$ or $p \equiv -1\pmod{4},$ then
$Z_p(M,\c)$  is a periodic map with period $2ps.$
If  $p \equiv 2\pmod{4}$ or $p \equiv 1\pmod{4},$ then
$Z_p(M,\c)$ is a periodic map with period $4ps.$\endproclaim

\demo{Proof} Let $\Si$ be the fiber.
Let $E$ be the associated fundamental domain of the infinite cyclic cover of
$M$
and $E_s = \cup_{0\le i \le s-1}T^i E$ as in \S 1. $E_s$ is diffeomorphic to
$\Si\times[0,1],$ forgetting \ps.
Thus $(Z_p(E))^s=Z_p(E_s)$ is then given by a scalar multiple, say $c$, of the
identity. By (3.5), $M_s$  has an induced \ps with $\si(\a(M_s)) = - 3\
\text{def\ }(M,\c_s).$  $<M_s>_p= c \dim V_p(\Si),$  the trace of $Z_p(E_s).$
Whereas if we gave $\Si \times[0,1]$ the product $p_1$-structure, then
$Z_p(\Si\times[0,1])$ would be the identity, and by (3.9), the associated
mapping torus, with a \ps with $\si$ equal to zero, would have $<\ >_p$ equal
to     $\dim V_p(\Si).$  Thus  $c=\k^{-9\text{def\ }(M,\c_s)}$  is a power of
$\k^3.$ If $p \equiv 0\pmod{4}$ or $p \equiv -1\pmod{4},$ this is a $2p$th root
of unity. If $p \equiv 1\pmod{4}$ or $p \equiv 2\pmod{4},$ this is a $4p$th
root of unity. So in the first case  $(Z_p(E))^{2ps}$ is the
identity. In the second case $(Z_p(E))^{4ps}$ is the identity. \qed\enddemo

\head \S 4  Links in $S^1 \times S^2$ and their
wrapping numbers  \endhead

Now we consider $\hat Z_{ p}(M,\c)$ where $M= S^1 \times S^2$
containing a  banded link $L$, and $\c$ evaluates to one on the
$S^1$ factor. We let $\hat z_{p}(L)$ denote this invariant. We will also let
$z_{p}(L)$ denote $Z_{ p}(M,\c).$ It turns out that for almost all $p,$  $\hat
z_{p}(L)$
is actually the reduction of a single automorphism of a free
$\Bbb Q[A,A^{-1}]$-module.  To see this we first  define a finite linearization
over
 $\Bbb Z[A,A^{-1}]$ of a certain weak cobordism category $\eusm C.$ By a weak
cobordism category, we mean a cobordism category which does not have a disjoint
union operation. The involution on the ring $\Bbb Q[A,A^{-1}]$ and
 $\Bbb Z[A,A^{-1}]$ sends $A$ to $A^{-1}$ and fixes $\Bbb Q$.

We form $\eusm C$ by taking the definition of the category $C^{p_1}_2,$
throwing out any mention of $p_1$ structure, insisting that every object be
diffeomorphic to either $S^2$ with an even number of banded points or the
$\emptyset,$ and insisting that every morphism be diffeomorphic to
either $S^3,$  $D^3,$  $S^2 \times I$ or $\emptyset$ equipped with a banded
link which meets each boundary component in an even number of points.

Now the Kauffman bracket on banded links in $S^3,$ is an involutory $\Bbb
Z[A,A^{-1}]$-valued invariant of closed objects of $\eusm C$ and so defines a
linearization $(\eusm V, \eusm Z)$ over $\Bbb Z[A,A^{-1}].$
If $\Si$ is a nonempty object with $2 n$ banded points   then $\eusm V(\Si)$ is
 the Kauffman skein module of $(B,2n),$ where $B$ is a 3-ball with boundary
$\Si.$
This module may be identified with
 $K_n,$  the Kauffman skein module for the disk with
$2n$-boundary  components discussed by Lickorish \cite{L2}.
This has a basis $\{D_i\}$  consisting of all the isotopy classes of
configurations of $n$ arcs in  $D^2$ with boundary the
collection of $2n$ points with diagrams with no  crossings.
The number of such diagrams and thus the dimension of this
$\Bbb Z[A,A^{-1}]$-module is the $n$th Catalan number $c(n)=\frac {1}{n+1}
\binom {2n}n$.
  One also has $\eusm V(\emptyset)=\Bbb Z[A,A^{-1}].$  Thus $\eusm V(\Si)$ will
be free and has finite rank.

If $M= S^1 \times S^2$ containing an even banded link $L$, and $\c$ evaluates
to one on the $S^1$ factor, then in the construction of \S1 we may always take
$\Si$ to be an object of $\eusm C,$ and $E(\Si)$ to be a morphism of $\eusm C.$
In this way we obtain an endomorphism $\eusm Z(E)$ of $\eusm V(\Si).$
The results (1.1)-(1.5), (2.3)-(2.6) apply just as well to finite linearization
of a
weak cobordism category.

Given a link $L$ in $S^1 \times S^2$, we may isotope it so
that it lies in  $S^1 \times B^2$, is transverse to $\{ 1\}\times
S^2$ and has a regular  projection to $S^1 \times B^1$. Next we
cut the diagram of the projection to  $S^1 \times B^1$ along
$\{ 1\}\times B^1$, to obtain a link diagram $\Cal  T$ in  $I
\times B^1$.  By the number of strands of $\Cal T,$ we  mean
the  number of points in $\Cal T \cap \{0\} \times B^1$.
Suppose for now that  $\Cal T$ has $2n$ strands.  Let $\Cal
D_i$ be a diagram for  $D_i$ in $[- 1,0] \times B^1$ with the
$2n$ points on $\{0\} \times B^1$.  Let $\Cal Q(\Cal T)$ be the
$c(n) \times c(n)$ matrix over $\Bbb Z[A,A^-1]$,  whose
$(i,j)$ entry is given by the coefficient of the skein element $D_j$
when  $D_i\cup \Cal T$ in  $[-1,1] \times B^1$ is written in terms of the basis
$\{D_j\}$  with the
$2n$ points now  on $\{1\} \times B^1$.   Let $\G(\Cal T) \in
\Bbb Z [A, A^{-1}]$  denote
the normalized characteristic polynomial of $\Cal Q(\Cal T).$

There is an alternative to the above method method for writing out $\Cal Q(\Cal
T).$ Let $D(n)$ denote the   $c(n) \times c(n)$
matrix whose $(i,j)$ entry is the bracket of the diagram in
$S^2$ obtained  by taking the union of the diagram for $D_i$
with the diagram for $D_j$   along their boundary. As this is a
diagram  without crossings this entry is  just
$\d=-(A^2+A^{-2})$ raised to the number of components in
the  resulting diagram in $S^2$. This matrix was first
considered by Lickorish  \cite{L2}.   $\det D (n)$ is a nonzero polynomial in
$\d$. It is easy to see that the diagonal entries of $D(n)$ are just $\d^n$ and
the off diagonal entries are $\d$ to smaller powers. Thus the degree of
$\det D (n)$ in $\d$ is $n c(n)= \binom {2n}{n}$. Thus $D(n)$ is invertible
over
 the field of rational functions $\Bbb Q (A).$ Let $B(\Cal T)$ be the matrix
over $\Bbb Z[A]$,  whose entries are given by the bracket polynomial of the
diagram $D_i\cup \Cal T\cup m(D_i)$ in $[-1,2] \times B^1$. Here $m(D_i)$ is
the diagram in $[1,2] \times B^1$ obtained by reflecting $D_i$ across the line
$\{1/2\} \times B^1.$  Note if $\Cal T $ is a diagram consisting of
$2n$ straight strands, then $B(\Cal T)= D(n).$ In general, we have
$\Cal Q(\Cal T)= B(\Cal T)D(n)^{-1}.$ Note that this method of calculating
$\Cal Q(\Cal T)$ does not make  it apparent that the entries of $\Cal Q(\Cal
T)$ lie in
$\Bbb Z [A, A^{-1}].$

We now consider the quantization $(\hat \eusm V,\hat \eusm Z)$ over the
rational functions $\Bbb Q (A).$ Let $\Cal Q(\Cal T)_\flat$ denote the induced
automorphism of $(K_n \otimes \Bbb Q (A))_\flat,$ as in \S1.  $\G(\Cal T)$ is
the characteristic polynomial of $\Cal Q(\Cal T)_\flat.$  Thus we have:

\proclaim{Theorem (4.1)} The similarity class of $\Cal Q(\Cal T)_\flat,$
and the polynomial $\G(\Cal T)$ are invariants of $L.$ \endproclaim

Thus we may let $\G(L)$ denote $\G(\Cal T),$ and let $D(L)$ denote the constant
term of $\G(\Cal T).$
The wrapping number $w(L)$ is the minimum number of transverse intersections of
$L$ with an essential embedded 2-sphere \cite{Li}.

\proclaim{Corollary (4.2)}If L is an even link in $S^1 \times
S^2$ with  diagram $\Cal T$ with $2n$ strands where  $n\ge
2$ and $\det (B(\Cal  T))$ is nonzero, then $w(L)=2n.$
\endproclaim

\demo{Proofs}  If $\det (B(\Cal T))$ is nonzero, then $\det
(B(\Cal  T)D(n)^{-1})$ is nonzero. Thus $\G(L)$ has degree
$c(n)$. Since $m<n$  implies $c(m) < c(n)$ for $n>2$, the
conclusion follows. \qed \enddemo

\proclaim{Corollary (4.3)} If L is an even link in $S^1 \times
S^2$, then  $c(\frac{w(L)}{2}) \ge \text{deg} (\Cal G(L)).$
\endproclaim

\demo{Proof} We may calculate  $\Cal G(L)$ from a tangle with
$w(L)$ of strands.\qed \enddemo

Hoste and Przytycki have calculated the Kauffman skein
module of $S^1  \times S^2$ and have in this way obtained results
on the wrapping number   \cite{HP}. Our results appear to be
different but a detailed comparison has not been  done. Hoste and Przytycki
show the Kauffman skein
module modulo $\Bbb Z[A,A^{-1}]$-torison is $\Bbb Z[A,A^{-1}],$
and let $\p$ denote the quotient map. In \cite{Gi}, we show that the trace of
$\Cal Q(\Cal T)$ is the same as $\p(L).$

Now we consider  an intermediate  linearization $(\check \eusm V,\check  \eusm
Z)$ over the PID $\Bbb Q [A,A^{-1}].$  Let $\Cal Q(\Cal T)_\sharp$ denote the
induced 1-1 endomorphism of $(K_n \otimes \Bbb Q [A,A^{-1}])_\sharp,$ as in
\S2.
$\G(\Cal T)$ is its characteristic polynomial.
Thus $\Cal I_{\check Z}(M,\c)$ is the principle ideal generated by $D(L).$
Thus in this case $k_{\Cal I}$  is $\Bbb Q [A,A^{-1},\frac{1}{D(L)}].$
Let $\Cal Q(\Cal T)_\natural$ denote the induced automorphism of $K_n
\otimes \Bbb Q [A,A^{-1},\frac{1}{D(L)}]$ as in \S2. We have :

\proclaim{Theorem (4.4)} The similarity class of $\Cal Q(\Cal T)_\natural,$
 over $\Bbb Q [A,A^{-1},\frac{1}{D(L)}]$ is an invariant of $L.$ \endproclaim

So we may let  $\eurm z(L)$ denote the similarity class of
$\Cal Q(\Cal T)_\natural,$  over $\Bbb Q [A,A^{-1},\frac{1}{D(L)}].$

$$\epsffile{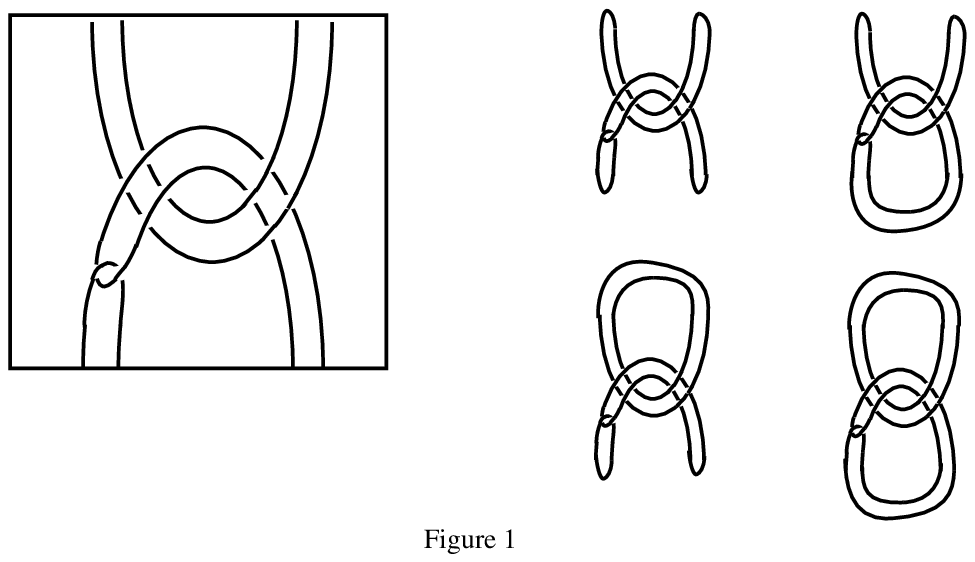}$$

\noindent {\bf Example (4.5)} Let $\Cal T$ be the tangle on the
left of  Figure 1. It can be closed up to form a link $L$ in
$S^1 \times  S^2$. $B(\Cal T)$ is given by taking the bracket
polynomial of the matrix  of link diagrams on the right of
Figure 1. Thus  $B(\Cal T)=\bmatrix  \d h & \d^2 A^6 \\
\d^2 h & w  \endbmatrix  $, where $h = -\d(A^4+
A^{-4})$, the  bracket of the standard diagram for the Hopf
link, and  $$w= -2 + A^{-16} - A^{-8} - A^{-4} - 2A^4 - 2A^8
- A^{20}.$$ Kauffman has a good method  for
calculating the bracket of link  diagrams with double strands \cite{KL,(\S
4.4)}. See also \cite{K2}, an early version of \cite{K3}.
We call this the Kauffman double bracket  method.  We used this method to
calculate
$w.$ Thus $\Cal Q(\Cal T)$ is given by
$$\bmatrix -1 - A^{-4}&-A^{-2} + A^6\\
 A^{-10} - A^{-6} + 3A^{-2} + A^2 - A^6 + 2 A^{10} - A^{14}&
  \ \ A^{-12} - A^{-8} +2 - 2 A^4 + A^{12} - A^{16} \endbmatrix  $$
This matrix has a nonzero determinant
$$D(L)= - A^{-16} + A^{-12} +2 - 2
A^4 -  A^{16} + A^{20}.$$
 Thus the above matrix represents $\eurm z(L)$, and
  $\G(L)
=  g_0 + g_1 x+x^2$ where $g_0 = D(L)$  above and   $$g_1=  -A^{-12} +
A^{-8} + A^{-4} -1 + 2 A^4 - A^{12} + A^{16}.$$
We may conclude that the wrapping number of $L$ is
four. This also follows from \cite{Li}, as well as \cite{HP}.  Using (1.5) and
an analog of (3.4), we conclude that $L$ is not isotopic
to its  image under a orientation reversing diffeomorphism of $S^1 \times S^2$
ignoring banding.

 If $X$ is either a scalar in $\Bbb Z[A]$ or a matrix  over $\Bbb Z[A],$ then
we let either $_pX$ or $X_p$ (depending on where there is more room for the
subscript) denote $X$ after evaluating at $A=A_p$. Similarly if $X$ is a
module or module homomorphism over
$\Bbb Z[A, A^-1]$ or $\Bbb Q[A, A^{-1}]$ let $X_p$ denote the result of
tensor product with $k_p$ or Id$_{k_p}.$

 We will say $p$ {\it is ordinary with respect to $n$} if and only if $D(n)_p$
is nonsingular. If $p$ is not ordinary, we will say it is {\it special with
respect to $n$}.  Ko  and Smolinsky \cite{KS} studied the question: when is $p$
ordinary with respect to $n?$      We note first that since $\det D (n)$ is a
nonconstant polynomial almost all $p$ are ordinary with respect to a given $n$.

Ko  and Smolinsky  showed that all the roots of $D(n)_p$
are of the form $\d= 2 \cos( \frac{k \p}{m+1})$ where $1 \le k \le m \le n.$
We note that one and two are ordinary with respect to any $n$.   Ko  and
Smolinsky showed that $2r$ is special with respect to $r-1,$ as required by
Lickorish \cite{L2}.

By \cite{BHMV1,(1.9)},  there is an epimorphism
$\e_{(D^3,2n)}: K(D^3,2n)_p\allowmathbreak \rightarrow
V_p(S^2,2n).$  Here we let $(S^2,m)$ denote the 2-sphere with
$m$ framed points. Using the nonsingular Hermitian form on $V_p(S^2,2n),$
 one sees that $\e_{(D^3,2n)}$ is an isomorphism if
and only if  $p$ is ordinary with respect to $n$.  In this case
$\{D_i\}$ describes a  basis for $V_p(S^2,2n)$ with cardinality  $c(n).$
\cite{BHMV1,(4.11 \& 4.14)} gives bases for
$V_p(S^2,2n)$ and so may also be used to determine
when $\e_{(D^3,2n)}$ is an isomorphism.   We obtain:

\proclaim{Proposition (4.6)} If $p\ge 4$ is even, then  $p$ is ordinary with
respect to $n$ if and only if $p > 2n+2.$ If $p\ge 3$ is odd, then  $p$ is
ordinary with respect to $n$ if and only if $p > n+1.$  \endproclaim

 \proclaim{Theorem (4.7)}
 If $p$ is ordinary with respect to $n$, and
$D(L)_{p}\neq 0$, then $\G_p(L)
=(\G(L))_p,$ and $\hat z_{p}(L)=\eurm z(L)_p$   \endproclaim

\demo{Proof}   If $p$ is ordinary with respect to  $n$,
$\hat Z_p (I \times S^2,\Cal T)$  with respect to the basis $\{D_i\}$ is
represented by the matrix $\Cal Q(\Cal T)_p.$  If
$D(L)_{p}\neq 0$,
then
$(\Cal Q(\Cal T)_\natural)_p$ represents $\hat Z_p (I \times
S^2,\Cal T)_\flat.$   \qed \enddemo

We note that the hypothesis of (4.7)  is true for almost all $p$ if $\G(L) \ne
0.$
Although it would  be interesting to calculate $z_{p}(L),$ or
$\hat z_{p}(L)$ for $p$  which are special with respect to $n$,
we do not pursue this now. In order to define invariants for odd links in $S^1
\times S^2$ we have several options. One could color the link $L$ with a fixed
even integer, say two, and evaluate the above invariants in the colored
theories. Alternatively one could consider the link obtained by replacing each
component of $L$ with two strands using the banding to form a new banded link
$L'$, and then calculate the invariants of the even link $L'.$ Note that  $L'$
is formed by taking a pair of ``scissors'' and splitting each band in $L$ to
form two bands. One could use other ``even'' satellite constructions.
In \S 10 we will mention a third method of obtaining  invariants of odd links.

\head \S 5 Knot invariants \endhead

Given an oriented knot $K$ in a homology sphere $S$, we may
let $S(K)$  denote zero framed surgery to $S$ along $K$.
$S(K)$ has the integral  homology of $S^1 \times S^2$.  We let
$\c$ denote the cohomology class  which evaluates to be one
on a positive meridian of $K$.      Let $Z_p (K)$ denote $Z_p (
S(K),\c),$ $\hat Z_p (K)$ denote $\hat Z_p (
S(K),\c),$   and  $\G_p(K)=\G_p (S(K),\c)$.  We will also let
$\Cal E(K)$ denote the list of eigenvalues of $\hat Z_{p} (K)$ counted
with  multiplicity. By Proposition (1.5),
one can see that these invariants do
 not depend on the string orientation of the knot.
Let $- K$ denote the  knot obtained
by taking the  mirror image of $K$
and reversing the string orientation. This knot represents the inverse of $K$
in the knot cobordism group. Then
$Z_p (-K)= (Z_p ( K))^*$. Let $U$ denote the unknot in $S^3$,
then  $S^3(U) = S^1 \times S^2$.  So far for the examples we
have calculated $\hat Z_{p}(K)$ is diagonalizable. It would be interesting to
find
two knots $K_1$ and $K_2$  such that
$\G_p(K_1)=\G_p(K_2),$ but $\hat Z_p(K_1)\ne \hat Z_p(K_2),$ or
$Z_p(K_1)\ne Z_p(K_2).$

\proclaim {Theorem (5.1)} If K is a fibered knot in a homology
sphere  which is a homotopy ribbon knot, then one is a root of
$\G_p (K).$  \endproclaim

\demo{Proof} According to Casson and Gordon \cite{CG1}, a
fibered knot  in a homology sphere is homotopy ribbon if and
only if the the closed  monodromy extends over a
handlebody $H$. In this case the ordinary  mapping torus $R$
of the extension to the handlebody is a homology $S^1 \times
B^3$ which embeds naturally in a homology ball $B$ with
boundary $S$ in which $K$ bounds a homotopy ribbon disk
$\D$. In fact $R$ is basically the exterior of $\D$.  We can  give
$B$ a $p_1$-structure..  It will  induce a \ps on $R$ which in turn induces a
\ps on  $S(K)$ with $\si$ zero. Note $H$ represents an element
in $V(\partial H)=  V(\Si)$ which is fixed by $T$. Thus $H$
represents an eigenvector with  eigenvalue one. \qed
\enddemo

\noindent {\bf Remark} Although there is an algorithm to
answer the question  of whether a given diffeomorphism
extends over a handlebody \cite{CL},  consideration of the
eigenvalues of a map induced  under a TQFT functor    may
turn out to be a good way to show that a diffeomorphism does
not  extend over a handlebody or even bound in the bordism
group of  diffeomorphisms \cite{B},\cite{EW}.

The following proposition  follows instantly from
\cite{BHMV1,(1.5)}. Here $i_p: k_2 \rightarrow k_{2p}$ and
$j_p: k_p\rightarrow k_{2p}$ are  the
homomorphisms defined in \cite{BHMV1}. We must note  that
for $p$ odd, $i_p(\k_2)j_p(\k_p)= \k_{2p}$. We use the same
symbols  to describe the induced maps on similarity classes
and polynomials over  these rings.
We will discuss the tensor product of polynomials in the Appendix.

\proclaim{Proposition (5.2)} If $p$ is odd, then
$Z_{2p}(K)=i_p(Z_2(K))  \otimes j_p(Z_p(K)),$ and so
$\G_{2p}(K)=i_p(\G_2(K)) \otimes  j_p(\G_p(K)).$ \endproclaim

Let $U$ denote the unknot. We will say the similarity class of the identity on
a free module of rank one is trivial.

\proclaim{Theorem (5.3)} For  all $p,$ $Z_p(U)$ is trivial.  If $p$ is one,
three or four, and $K$ is a
knot in $S^3$, then  $Z_p(U)$ is also trivial. \endproclaim

\demo{Proof}  The first statement follows directly from the definitions.
 If $p$ is one,
three or four, we have \cite{BHMV1,\S 2}\cite{BHMV2,\S 6} $\o=\eta$,
 $\k^{-3} \eta=1,$  $\k^6=1,$ and thus $\k^ 3\eta=1.$ One may
obtain any  knot in $S^3$ from the unknot $U$ by doing $\pm
1$ surgery to the  components of an unlink in the
complement of $U$ where the linking  number of each
component with $U$ is trivial \cite{R,(6D)}.  Thus we may pick
a Seifert surface $F$ for $U$ in the complement of this
unlink.    We may calculate $Z_p(U)$ from an $E$ that we
construct with  $\Si$ = $F$ capped off. Here we give $S^3$ a \ps
with $\si$ zero.  Let us  calculate the effect of a single surgery
on $E$.  Let $E'$ denote the result  of performing $+1$ framed
$p_1$-surgery to $E$. Then $Z_p(E') = \o_p  Z_p(E)=\eta
Z_p(E)$.  Let $K'$ denote the image of $U$ after this  surgery.
Then  $Z_p(K') = \k^{-3} Z_p(E')=Z_p(E)=Z_p(U).$ If we  were to
perform $-1$ surgery then the invariant would change by a
factor  of $\k^ 3\eta=1.$ This same argument may be repeated
to show the further  surgeries do not change $Z_p(K).$  \qed
\enddemo

In order to obtain (5.6) below about periodicity of $Z_p(K),$ for $K$ fibered
with periodic monodromy, we make the following definitions and observations
which will be useful later as well.
Let  $\sigma_{\omega}(K)= \text{Sign} ( (1- \o) V+ (1- \bar\o) V^t),$
where $\o \in \Bbb C$ with $|\o|=1$ and $V$ is a Seifert matrix for $K.$
Following \cite{KM2}, let  $ \si_d(K)= \sum_{i=1}^{d-1} \sigma_{\o_d^i}(K),$
where $\o_d=e^{2 \p i/d}.$
These are called  the total $d$-signatures of $K$. Our convention is that
$\si_1(K)$ is taken to be zero.  Also $\text{def}(S^3(K),\c_d)=-\si_d(K).$  By
(3.5) we have:

\proclaim{Proposition (5.4)} Let $\a(K,d)$ be the \ps on $S(K)_d$ induced from
$S(K),$ $\si(\a(K,d))= 3 \si_d(K).$ \endproclaim

\proclaim{Lemma(5.5)} If $K$ is a fibered knot with a periodic monodromy of
order $s$, then $\si_{s+1}(K)\equiv 0 \pmod {8},$ and $\si_s(K)\equiv 0 \pmod
{8}.$
\endproclaim
\demo{Proof} Let $D_d$ denote the $d$-fold cyclic branched cover of $D^4$ along
a pushed in Seifert surface for $K$. $D_d$ is a spin simply connected manifold
with boundary $K_d$ and Sign$(D_d)= \si_d(K).$
 $K_{s+1}$ is  the result of $1/n$ surgery on $K$ for some $n.$
Thus $K_{s+1}$ is a homology sphere.   So the signature of $D_{s+1}$ is the
signature of an even  unimodular symmetric matrix and so is zero modulo eight.
$(S^3(K))_s$ is diffeomorphic to $\Si \times S^1.$ It follows that $H_1(K_{s})$
is torsion free. Thus the induced form on $H_2(D_s)$
modulo the radical of the intersection form is given by an even  unimodular
symmetric matrix and so has signature zero modulo eight.
\qed\enddemo

\proclaim{Proposition (5.6)}  Suppose
 $K$ is a fibered knot with a periodic monodromy of order $s$. Let
$u^{4(\si_s(K)/8)}$ have order $h$ in $k_p.$ Note $h$ divides $p.$
Then $Z_p(K)$ is a periodic map with period $hs.$\endproclaim

\demo{Proof} By(5.4), $\k^{\si(\a(K,d))}=u^{4(\si_s(K)/8)}.$ The result follows
from the proof of (3.10).\qed\enddemo

$$\epsffile{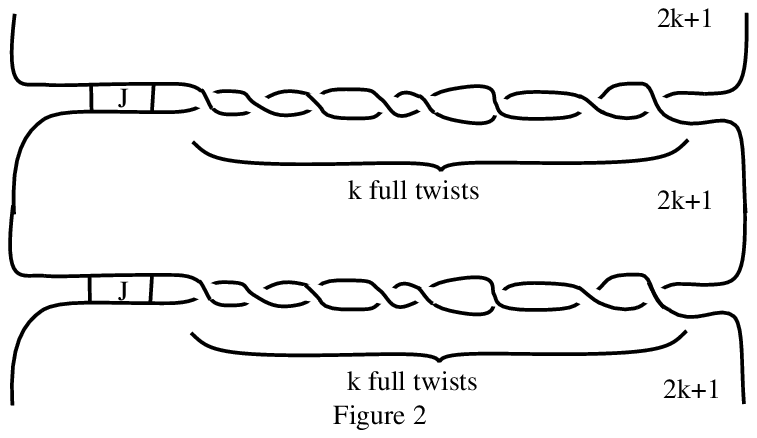}$$

Given a surgery description of a knot as in Rolfsen
\cite{R,159} where  each surgery curve has zero linking
number with the unknot, we will give  a procedure to
calculate $Z_p (K)$. We use Rolfsen's method of giving a
surgery description for  the infinite cyclic branched cover of
$S^3$ along  $K$ from a surgery description for $K$ \cite{R,
p162, p158}. The same  picture and argument shows that we
are obtaining a description of the  infinite cyclic unbranched
cover of $M(K)$ as obtained by surgery to  $\Bbb R \times
S^2$.  Since we are actually working with manifolds with a
\ps our surgery is $p_1$-surgery, but when we do our initial
surgery to tie  up $K$ we change the \ps on $S^3$ and $S^1
\times S^2$. In other words  our surgery description
describes $M(K)$ with the \ps with $\frac{\si}3$  equal to the
number of plus one surgeries minus the number of minus
one  surgeries done \cite{BHMV1, Appendix II}.  We have
$Z_p(K) = \k_p^{- \si} Z_P(E'(\Si))$, where $E'(\Si)$ is the
manifold we describe below.

\subhead Twisted Doubles of Knots
\endsubhead

We consider first the case that $K$ is $D_k(J)$, the $k$-twisted
double of  a knot $J$ \cite{R,p112}. We perform a single minus
one surgery to undo  the clasp.  See \cite{CG2}, \cite{K1,Chapt
18} where the (finite) cyclic  covers are calculated for $J$  the
unknot, and \cite{R} for $k=0$ and $J$  equal to the right
handed trefoil. Figure 2 gives a surgery description for  the
infinite cyclic cover (with $k=4$). The box labelled $J$
represents two  parallel copies of a string diagram for $J$ with
zero writhe. To obtain the  infinite cyclic cover one should
perform framed surgery to $\Bbb R \times S^2$ along the
indicated infinite chain. Consider a 2- sphere given by $\{
a\}\times S^2$ which meets one component of our  diagram in
two points. Delete these two points and close the surface off
by  adding a  tubular neighborhood of an arc on the
intersected component.  Call the resulting torus $\Si.$
$$\epsffile{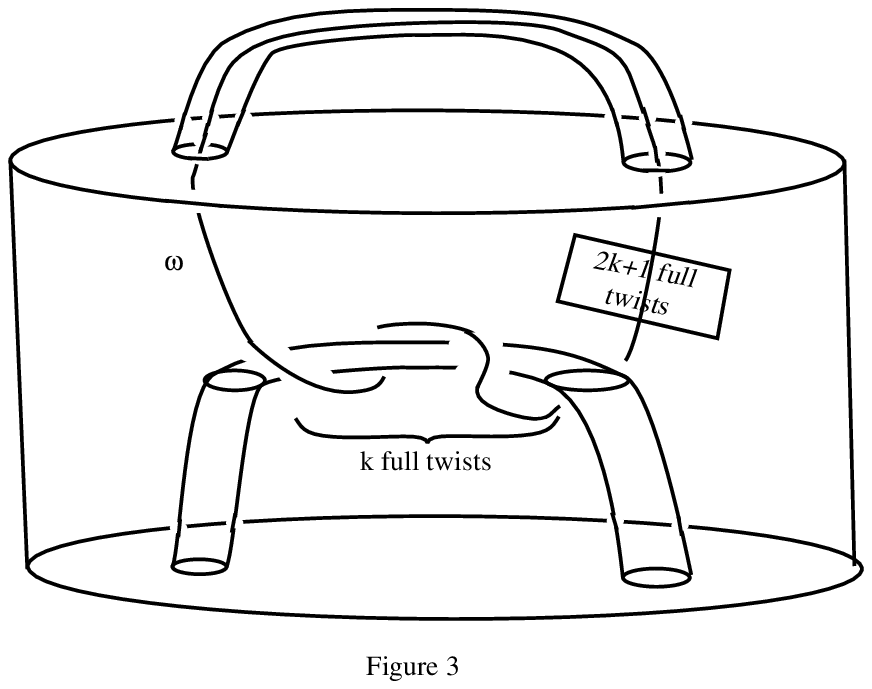}$$

We can now construct  $E'(\Si)$ by taking a  ``slab'' $I \times
S^2$,   drilling out  a tunnel along an arc which meets $\{0\}
\times S^2$ (the arc  is the top of one of the  surgery
curves), adding a 1-handle added  along  $\{1\} \times S^2$,
and finally performing framed surgery along  a simple closed
curve $\g$ which travels once over the one handle.
According to \cite{BHMV1,(1.C),\S 2,(5.8)}, $Z_p(E'(\Si)) =
Z_p(E''(\Si))$, where $E''(\Si)$ is formed by   placing a linear
combination of banded links $\o$ with $2k+1$ full twists along
$\g$ where  we would have performed surgery in
constructing $E'(\Si).$  In Figure 3,  we sketch $E''(\Si)$ when
$J$ is the unknot. The figure actually shows the  case $k=2.$
We only draw an $I \times B^2$ portion of  $I \times S^2.$
Given an element in $\a \in V(\Si),$  we may figure out its
image under  $Z_p(E''(\Si))$ as follows.  $\a$ can be
represented by a banded link in a  solid genus one
handlebody $H$ with boundary $\Si$. If we glue this
handlebody to the bottom of $E''(\Si)$ we will obtain a linear
combination  of banded links in the handlebody $H \cup
E''(\Si)$ with boundary  $T(\Si)$. This linear combination is
formed as the union of $\a$ and $\o$.  This linear combination
represents $Z_p(E''(\Si))(\a) \in V(T(\Si)).$

\subsubhead Preparatory material on the Kauffman bracket and the Kauffman
polynomial
\endsubsubhead

Let $c_k(J)$ denote the bracket of the  diagram, say $\Cal D_k(J),$ obtained
from a diagram
of a knot
$J$ with writhe $k$ by replacing each arc by 2 parallel arcs.
The Kauffman double bracket
method is a very efficient means for  calculating  $c_k(J).$  In fact
$c_k(J)= [\Cal D_k(J)]_2+1,$ \cite{KL,(p.35-36)} and $[ \ \ ]_2$ satisfies a
simple skein relation.  We let
$[[J]]$ denote $[\Cal D_0(J)]_2.$
   Note $[[ J_1\#J_2]] = \frac{1}{\d^2-1} [[ J_1]]\ [[ J_2]].$  Also $c_k(J)=
A^{8k}[[J]]+1.$

If $\Cal D$ is a knot diagram, then the writhe of $D$, denoted $w(\Cal D)$
is well  defined. This is not true of link diagrams.  One may
show that $<\Cal D>_2 =  i^{w(\Cal D)} 2.$  If $J$ is a knot, let $<J>$ denote
the bracket polynomial of a diagram for $J$ with zero  writhe.
Thus $<J>_2 =  2.$ Using the relation of the bracket
polynomial   and the Jones polynomial, it is easy to see that
$<J>$ will be a polynomial in even powers of $A$. It is convenient to note that
\cite{K5,(3.2)}
$$<J>= \left(\frac{a+a^{-1}}{z}-1\right) F_J(a,z)
_{| a= -A^3\text{ and } z= A+A^{-1}}$$

Here $F_J(a,z)$ denotes the Kauffman polynomial, normalized so that it is one
for the unknot. This is the form that is in the tables of \cite{K1}.
As observed in \cite{KL},
$[\Cal D]_2$ is the Dubrovnik version of the Kauffman polynomial
\cite{K5,\S VII} of $\Cal D$  evaluated at $z= A^4 -A^{-4},$ and $a= A^8.$
Using Lickorish \cite{L4}, we have:
$$[[J]] = -\left(\frac{a+a^{-1}}{z}-1\right) F_J(a,z)
_{| a= - i A^8\text{ and } z= i(A^4 - A^{-4})}$$

Let $\eurm b_k(J)$ denote the bracket of $\Cal D_0(J)$  but with $k$ additional
full
twists between the two strands.  One has that $$\eurm b_k(J)= A^{-6k} c_k(J)
=A^{2k}[[J]]+A^{-6k}.$$
The following proposition
now follows  easily from the above formulas. This proposition
seems related to identities  in \cite{P} and could perhaps be
proved using them.

\proclaim {Proposition (5.7)}$_p\eurm b_k(J)$ is periodic in $k$ with
period  $p.$    Also $_2\eurm b_k(J)= (-1)^k 4.$  \endproclaim

\subsubhead Some conventions which hold for the rest of this
paper  \endsubsubhead

   We will let $A$ denote $A_p$, unless there may be some
confusion as to  which $p$ is meant.  We do not specify which
primitive $2p$th root of  unity this is.  Similarly
 $\O$, $\o$, $\k$, $\d$ will  denote the items defined in
\cite{BHMV1}, and  sometimes denoted $\O_p$, $\o_p$, $\k_p$,
$\d_p.$ We let $\mu= \mu(1)  = -A^3.$ We also give $\k^{-3} \eta$ the
name $\b.$   We also let $\eurm b_k(J)$ denote  $\eurm b_k(J)$ evaluated at
$A_p$, and $<J>$ denote $<J>_p$

\subsubhead The case $p=2$  \endsubsubhead

 If $\Si$ is a torus, then $V_2(\Si)$ has a basis consisting of
$1$ and $z$. We have that $\o  = \eta   \O $ where $\O  = 1+
\frac z2,$ and $\b  = \frac{1-A}{2}$. With respect  to the basis
$\{1,z\}$, $Z_2(E(\Si)) = \k^{-3} Z_2(E''(\Si))$ is given by $\b
\bmatrix 1 & <J>\\ \frac{\mu^{2k+1}<J>} {2 \d}
&\frac{\mu^{2k+1} \eurm b_k(J) }{2 \d}  \endbmatrix.$  Thus
   $Z_2(E(\Si))$ is given by $\frac {1-A }{2}   \bmatrix 1 & 2\\
(-1)^k \frac{A } {2} &A \endbmatrix$.

\proclaim{Proposition (5.8)} Let $K$ be the $k$ twisted double of
$J$. If $k$ is even, then $Z_2(K)$ is trivial.  If $k$ is odd,
$Z_2(K)= \frac {1-A }{2}   \bmatrix 1 & 2\\
- \frac{A } {2} &A \endbmatrix_{\natural},$ and $\G_2(K)$ is $x^2-x+1.$  So
$\Cal  E_2(K)$ is $\{1\}$ if $k$ is even and is
$\{A_3,\bar A_3\}$ if $k$ is  odd.\endproclaim

\subsubhead The case $p=5$  \endsubsubhead

If $\Si$ is a torus, then $V_2(\Si)$ has a basis consisting of $1$
and $z$. We have $\o  = \eta \O $.   $\O  =  1+ \d  z,$  and $\b =
\frac{3 -A  + 4 A^2 - 2A^3}{5}.$ With respect to the  basis
$\{1,z\}$, $Z_5(E(\Si)) = \k^{-3} Z_5(E''(\Si))$ is given by $\b
\bmatrix 1 & <J>\\ \mu^{2k+1} <J> &\mu^{2k+1}   {\eurm b_k(J) }
\endbmatrix.$

By Proposition (5.7), each
entry in the above  matrix  has period five. Thus:

\proclaim {Proposition (5.9)} $Z_5(D_k(J))= \left(\beta \bmatrix 1 & <J>\\
\mu^{2k+1} <J> &\mu^{2k+1}   {\eurm b_k(J) }
\endbmatrix \right)_\natural$.  $Z_5(D_k(J))$ is periodic in $k$ with
period  five.
\endproclaim

We now calculate these five invariants for various knots $J$.
First we  consider $J =U$, the unknot. We can make a number
of predictions  a priori. Note that $D_0(U)$ is the unknot again,
so $Z_{5}(D_{5n}(U))$ is trivial.  Also note  $D_1(U)$ is the figure
eight knot and $D_{-1}(U)$ is right  handed trefoil. Both of
these knots are fibered, so $\Cal  E_{5}(D_{5n+1}(U))$  and
$\Cal E_5(D_{5n-1}(U))$ should consist of  elements of norm
one. Since the trefoil is period with period six, by (5.4), $\Cal
E_{5}(D_{5n+1}(U))$ should consist of 30th roots of unity.    Also since the
figure eight is amphichiral,  $\overline{\Cal
E_{5}(D_{5n+1}(U))}= \Cal E_{5}(D_{5n+1}(U)).$  We obtain:

\proclaim {Proposition(5.10)} If $k \equiv 0\pmod{5},$
$Z_{5}(D_{k}(U))$ is trivial. If $k\not \equiv 0\pmod{5},$
$Z_{5}(D_{k}(U))= \beta \bmatrix 1 & \d\\ \mu^{2k+1} \d &\mu^{2k+1}
(A^{2k}(\d^2-1)+A^{-6k})
\endbmatrix_{\natural}.$  In particular for all   integers n, we have:
  $$\eightpoint
\aligned
&\G_{5}(D_{5n}(U))= x-1\\ &\G_{5}(D_{5n+1}(U))= x^2- (A+
\bar A)x + 1\\ &\G_{5}(D_{5n+2}(U))= x^2- (1+ \bar A)x +
\bar A\\ &\G_{5}(D_{5n+3}(U))= x^2- (1 +\bar A^2)x + \bar
A\\ &\G_{5}(D_{5n+4}(U))= x^2- ( \bar A)x + (\bar A)^2\endaligned \qquad
\aligned &\Cal E _{5}(D_{5n}(U))= \{1\}\\ &\Cal E _{5}(D_{5n+1}(U))= \{
A, \bar A \}\\ &\Cal E _{5}(D_{5n+2}(U))= \{ 1, \bar A \}\\
&\ \\
&\Cal E _{5}(D_{5n+4}(U))= \{ A_3 \bar A, \bar A_3\bar A \}.
\endaligned$$
\endproclaim

Note $\Cal E _{5}(D_{5n+4}(U))$ are primitive 15th roots of
unity. We  do not list the eigenvalues for the $5n+3$ twisted
doubles. Although they  are easily worked out, the formulas
are not enlightening. By (5.1) and (5.10), we have that the trefoil and figure
eight knots are not homotopy ribbon. This a (not very deep) four dimensional
result obtained from studying a TQFT in dimension 2+1. By (7.6), below we have
that the granny knot is not homotopy ribbon as well.

Let $RT,$ $LT$ and  $F8$ denote the right handed trefoil, the left handed
trefoil and the figure eight knots.
For $J$ = $RT$, $LT$  $F8,$ and the square knot $RT \# LT$ and for all k, the
above matrix has nonzero determinant.
Below we list $\{k,\G_{5}(D_{5n+k}(J))\}$ for $J = RT,$ $LT,$  $F8,$ and $RT \#
LT.$

$$\eightpoint
\allowdisplaybreaks\align &J= RT\\
&\{ 0,1 + 2\,{A^2} - 2\,{A^3} -
 \left( 2 - A + 2\,{A^2} - {A^3} \right) \,x + {x^2}\}
\\
&\{ 1,-{A^3} - \left( 2 - {A^3} \right) \,x + {x^2}\}
\\
&\{ 2,1 + A - {A^2} - \left( 1 - A + {A^2} - 2\,{A^3} \right) \,x + {x^2}\}
\\
&\{ 3,1 - 2\,A - {A^3} - \left( 1 - A \right) \,x +
 {x^2}\}
\\
&\{ 4, -A + A^2 + {A^3}   - A\,x + {x^2}\} \\
 &J= LT \\&\{ 0,2 - 2\,A - {A^3} + {A^3}\,x + {x^2}\}\\
&\{ 1,-1 - A + {A^2} - \left( 1 - A + {A^2} \right) \,x + {x^2}\} \\
&\{ 2,1 + 2\,{A^2} - {A^3} - \left( 1 + A \right) \,x + {x^2}\} \\
&\{ 3, A - A^2 - {A^3}   - \left( 2 - A + 2\,{A^2} - 2\,{A^3} \right) \,x +
{x^2}\}\\
&\{ 4,1 - \left( 2 - A - {A^3} \right) \,x + {x^2}\}  \\
 &J= F8\\
 &\{ 0,-3 + 2\,A - 2\,{A^2} + 3\,{A^3} - {A^2}\,x + {x^2}\} \\
&\{ 1,3 + 2\,{A^2} - {A^3} - \left( 2 + {A^2} \right) \,x + {x^2}\} \\
&\{ 2,1 - 2\,A - {A^3} - \left( 2 + {A^2} - 2\,{A^3} \right) \,x + {x^2}\} \\
&\{ 3,1 + 2\,{A^2} - {A^3} - \left( 2 - 2\,A + {A^2} - 2\,{A^3} \right) \,x +
{x^2}\} \\
&\{ 4,1 - 2\,A - 3\,{A^3} + {A^2}\,x + {x^2}\}\\
&J= RT \# LT\\
&\{ 0,-6 + 4 A - 4 {A^2} + 6 {A^3} -
   \left( A - {A^2} + 2 {A^3} \right)  x + {x^2}\} \\
&\{ 1,6 + A + {A^2} - \left( 1 + A + 2 {A^2} - {A^3} \right)  x + {x^2}\}\\
&\{ 2,1 - 5 A + {A^2} - 2 {A^3} -
    \left( 4 - 2 A + 2 {A^2} - 2 {A^3} \right)  x +
    {x^2}\}\\
&\{ 3,2 - A + 5 {A^2} - {A^3} -
    \left( 1 - {A^2} - 2 {A^3} \right)  x + {x^2}\} \\
&\{ 4,-A - {A^2} - 6 {A^3} -
    \left( -2 A + {A^2} - {A^3} \right)  x + {x^2}\}\endalign $$

\subsubhead The case $p=6$  \endsubsubhead

By  (5.2), $Z_{6}(K)=i_3(Z_2(K)) \otimes j_3(Z_3(K)).$ Note that
$i_3(A_2) = A_6^9 =-A_6^3.$ As
$Z_3(K)$ is  trivial, and $i_p$ fixes $\Bbb Z,$ we have
  by (5.8) that:

\proclaim{Proposition(5.11)}  $Z_{6}(K)=i_3(Z_2(K))).$ In particular, if $k$ is
even $Z_{6}(D_k(J))$ is trivial. If $k$ is odd,
$Z_6(D_k(J))= \frac {1+A^3 }{2}   \bmatrix 1 & 2\\
 \frac{A^3 } {2} &-A^3 \endbmatrix_{\natural},$  and $\G_6(K) = x^2-x+1.$  So
$\Cal E_6(D_k(J))$ is
$\{1\}$ if $k$ is even and is $\{A_3,\bar A_3\}$ if  $k$ is odd.
\endproclaim

\subsubhead The case $p= 10$  \endsubsubhead

By (5.2), we have $Z_{10}(K)=i_5(Z_2(K)) \otimes
j_5(Z_5(K)),$ and so  $\G_{10}(K)=i_5(\G_2(K)) \otimes
j_5(\G_5(K)).$ Thus one may work  out, for instance, $\G_{10}
(D_k(U)$. We note that
 $j_5(A_5) = A_{10} ^6.$ By (5.8) and (5.10),
$i_5(\G_2(D_{5n+4}(U))) = x^2-x+1,$ and
$j_5(\G_5(D_{5n+4}(U))) = x^2+ ( A_{10} ^4)x + A_{10} ^8.$

  Using A.3 from the appendix, we obtain: $$\G_{10} (D_{5n+4}(U))= x^4+
(A_{10}^4)x^3  -(A_{10}^2)x -A_{10}^6.\tag 5.12$$

\subsubhead The general case, $p\ge 3$  \endsubsubhead

  If $p\ge3$ then $\o = \eta \O$, and $\O = \sum_{s=0}^{n-1}<e_s>
e_s$  where   $e_s$ is an eigenvector for the twist map with
eigenvalue $\mu(s)= (-1)^s A^{s^2+2s}$, and $<e_s>= (-1)^s \frac{A^{2s+2}
-A^{-2s-2}}{A^{2}
-A^{-2}}$. The $\{e_i\}_{i=0}^{n- 1}$ where $n=\left[
\frac{p-1}2\right]$ form a basis for $V$ of a torus. We form two
$n\times n$ matrices $\eurm B(J,k)$ and $\eurm L(J)$.  $\eurm
B(J,k)_{i,j}= \b\sum_{s=0}^{n-1}<e_s> \mu(s)^{2k+1} b(J,k)_{i,s,j}$
where $b(J,k)_{i,s,j}$ is the bracket polynomial of the colored banded link in
Figure 4a.
$$\epsffile{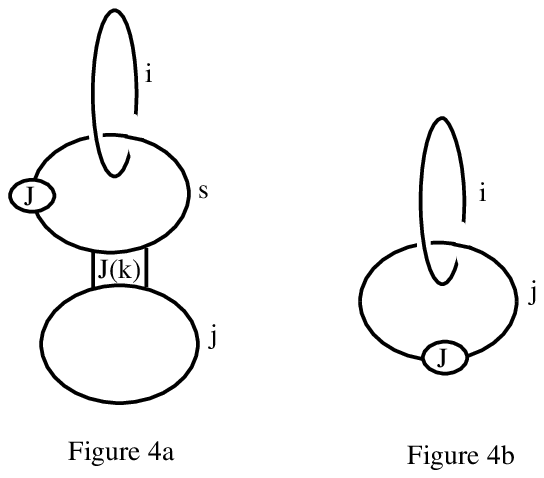}$$
The box labelled
$J(k)$ represents two parallel  strands of a string diagram for
$J$ with zero writhe with $k$ additional full twists added to the strands
and the circle labelled J represents a string diagram for
$J$ with zero writhe.  Let $\eurm L(J)_{i,j}$ be the bracket
polynomial of the colored banded link in
Figure 4b.
For a knot $J,$ let $J_c$ denote $J$ colored $c.$

\proclaim{Theorem (5.13)} If $p\ge 3$, and $<J_c>_p$ is nonzero for all $0\le
c\le n-1$,
then $Z_p(D_k(J))=(\eurm B(J,k)\  \eurm L(J)^{-1})_\natural.$\endproclaim

 $b(J,k)_{i,s,j}$ and  $\eurm L(J)_{i,j}$ may be
calculated recursively using the colored Kauffman relations. These are given in
\cite{MV2} and for p even in \cite{KL}. In fact $\eurm
 L(U)_{i,j}$ for $p$ even is given in \cite{KL,p.127}. This is easily worked
out for all p using the formulas in \cite{MV2}. It is not hard to see that the
summations given \cite{MV,p 367} should be taken over $k$ such that $(i,j,k)$
is a small admissible triple, when one specializes to $A=A_p.$  Using the
colored Kauffman relations, one may deduce:

\proclaim{Corollary (5.14)} If $p\ge 3$, and $<J_c>_p$ is nonzero  for all
$0\le c\le n-1,$
then $Z_p(D_k(J))$ is periodic in $k$ with period $p.$\endproclaim

Let $\d(k;i,j)$,
and $<i,j,k>$ be as in \cite{MV2}.   We have:

$$\eightpoint\align
\eurm L(U)_{i,j}&=  \sum_{
\text{$r$  $\ni$ $(i,r,j)$ is a small admissible triple} } \d(r;i,j)^2 <r>.\tag
5.15\\
b(U,k)_{i,s,j}&= \sum_
{
\matrix
\text{$r$,$r'$ $\ni$ $(i,r,s),$ and $(j,r',s)$}\\
 \text{are small admissible triples}
		\endmatrix
	}
\d(r;i,s)^2 \d(r';j,s)^{2k} \frac{<r> <r'>}{<s>}\tag 5.16 \endalign$$

Of course both of these should be evaluated  at
$A=A_p.$    We have calculated $Z_p(D_k(U))$ exactly with entries polynomials
in $A$ for $0 \le k \le p-1$ and for $5 \le p \le 16.$ These as well as other
lists of  quantum invariants are available at gopher://math.lsu.edu.
For certain knots, we have carried the calculation to higher $p.$  We observe
that for $p \le 20,$  $Z_p(RT)$ is a periodic map with period $3p,$  and
sometimes less. By (5.6), $Z_p(RT)$ must be periodic with period  by $6p,$ and
sometimes less,
since the trefoil has a monodromy of period six. The roots of
$Z_p(F8)$ are all periodic maps for $p\le 20$.
The  period is an  erratic function of $p$.  This periodicity is somewhat
surprising since the monodromy for F8 is hyperbolic. We also observe that
$\G_p(6_1)$ has one as an eigenvalue one for $p\le 18$. The tweenie knot  is
$5_2= -D_2(U)$  \cite{As,p.86}.  We noticed
that the degree of $\Gamma_{2r}(\text{ tweenie knot  })$ is less than $r-1 =
\dim V_{2r}(\text{ torus })$ for $3
\le r \le 9.$ Thus the hypothesis of (2.7) does not hold. We plan to
investigate whether $\Cal I(M,\chi)$ is principle in this case.

 \subhead Other Knots $K$  \endsubhead

If we consider more general knots, two extra difficulties arise.
First,  we  may not be able to calculate $Z_p(E(\Si))$ directly
but only some power of  it. As an example, starting with the surgery
description of the knot $8_{16}$ given in \cite{R,p.169}, one may calculate
 $Z_p(E(\Si))^n$ for $n \ge 2.$  Fortunately if we can calculate
$Z_p(E(\Si))^c$, then we can also calculate  $Z_p(E(\Si))^{c+1}$,
and from this deduce $Z_p(E(\Si))_{\natural}.$ Secondly, the
genus of $\Si$ may be higher than one. In the higher genus
case, there is a good description \cite{MV1},\cite{BHMV1} of a
basis of  $V_p(\Si)$ in terms of colorings of a trivalent graph  which is a
deformation retract of the  handlebody.   Thus the same
methods may be applied.  The calculation becomes more
difficult as $p$ and the genus of $\Si$ grow.
\subsubhead A computation with $\Si$ genus two \endsubsubhead

$V_5$ of the boundary of a genus two handlebody
has a basis given  by the links denoted $1,$ $z,$ $w,$ $zw,$ and $z\#w $ in
Figure 5. Starting with a surgery description of the knot $8_8,$ we obtain, as
in  Rolfsen, a surgery description of the infinite cyclic cover shown in Figure
6. Just as above, we can then find a matrix for calculating $Z_5(8_8).$ In
particular, we have that $\G_5(8_8)=(x-1) (x^4 + (-1 - A - A^2)x^3 +
 (1 + A^2)x^2 + (-A - A^2 - A^3)x+ (-1 + A + A^3)) .$

$$\epsffile{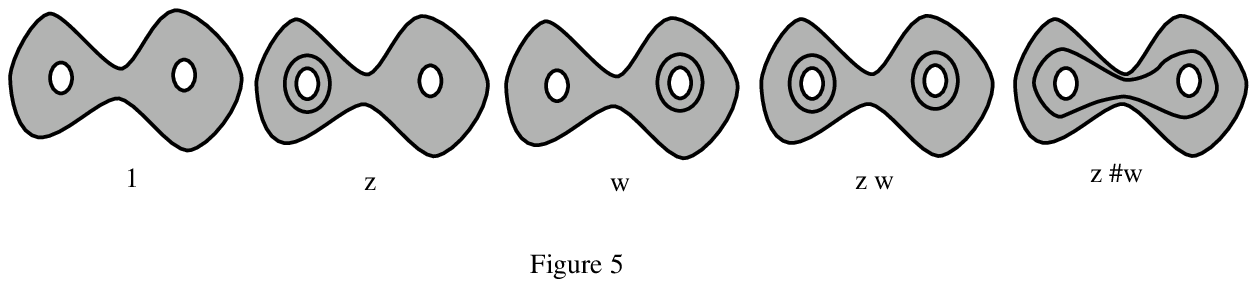}$$
$$\epsffile{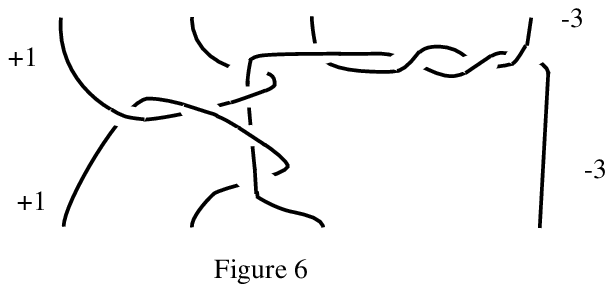}$$

\head \S 6  Quantum Invariants of the finite cyclic covers of $S^3(K)$
\endhead

For $p \ge 3$,  (1.8), and the remarks following tell us how  to compute
$<S^3(K)_d>_p$ recursively as a function of $d$, once we know $\G_p(K).$

We  discuss $<S^3(K)_d>_p$ for low $p.$
We have that $<S^3(K)_d>_p =1 \text{, for } p =1,3, \text{ and }
4.$ For $p=1$, this follows easily from the definitions.
For $p =3,$ and $4,$ this follows from (1.8) and (5.3).
By (1.7) and (5.9)
$$ <S^3(D_k(J))>_5 = \b_5 (1+ \mu_5^{2k+1} \eurm b_k(J)_5).\tag 6.1$$
By (1.8) and (5.11)
$$ <S^3(D_k(J))_d>_6 = \cases
1& \text{if $k$ is even}\\
2& \text{if $k$ is odd and $d \equiv 0\pmod{6}$}\\
1& \text{if $k$ is odd and $d \equiv \pm 1\pmod{6}$}\\
-1& \text{if $k$ is odd and $d \equiv \pm 2\pmod{6}$}\\
-2& \text{if $k$ is odd and $d \equiv 3\pmod{6}$}\endcases
				.\tag 6.2$$
We let $s(k,d)$ denote the function above.
Again as $Z_2$ does not satisfy the tensor product axiom, we need an
alternative
method of calculating $<S^3(K)_d>_2. $
 $<S^3(K)_d>_2$ can be computed from $<S^3(K)_d>_6$ using \cite{BHMV1,1.5}
since $<S^3(K)_d>_3 =1.$
In fact, one sees that $<S^3(K)_d>_2$ is equal to the sum of the $d$th powers
of the roots of
$\G_2(K).$ Thus it turns out that $<S^3(K)_d>_p$ is the sum of the dth powers
of the roots of $\G_p(K),$ for all $p\ge 1.$  In particular,
$$ <S^3(D_k(J))_d>_2 = s(k,d).\tag 6.3$$ Thus we have

$$ <S^3(D_k(J))_d>_{10} =  s(k,d) j_5(\b_5 (1+ \mu_5^{2k+1} \eurm
b_k(J)_5)).\tag
6.4$$

Here are some examples for $p=5.$  All of these examples are
atypical except perhaps the last. In the first of these examples, we compare
our result with previous calculations.

\subhead The covers of 0-surgery along the trefoil   \endsubhead

By (5.10) the eigenvalues of $Z_5(RT)$ are primitive 15th roots of unity. Thus
$<S^3(RT))_d>_5$ is periodic with period fifteen. Actually the first fifteen
values are $
-{A^4},
{A^3},
2\,{A^2},
A,
-1,
-2 A^{-1},
{A^3},
-{A^2},
-2\,A,
-1,
A^{-1},
-2\,{A^3},
-{A^2},
A,
2.$  Since the monodromy of the trefoil has order six, one
might, at first, expect periodicity of order six. However  by (5.4),
$<S^3(RT)>_5=\k^{-3\si_7(RT)}<(S^3(RT))_7>_5=
-A^4.$ This is consistent  with our previous calculation of $<S^3(LT)>_5.$
$S^3(RT)_6$ is the three torus. Thus
the invariant of the three torus equipped with a \ps $\a$ with $\si(\a)$ zero
is   $\k^{-3\si_6(RT)}<(S^3(RT))_6>_5= \k^{24}(-2A^{-1})=2.$
This calculation agrees with (3.9).

\subhead The covers of 0-surgery along the untwisted double of the figure eight
knot \endsubhead

$<S^3(K)_d>_5$ is given by the sum of the $d$th powers of the roots of
$\G_5(K)$ counted with multiplicity. It may also be easily computed
recursively. This example is  atypical in that we noticed something systematic.
$\G_{5} D_0(F8)=x^2- (A^2)x+(-3+2 A -2 A^2 +3A^3).$
Let $\l = 12 (A + \bar A) -8 (A^2 + \bar A^2)-1.$ $\l$ is a positive
real number under all complex embeddings of $\k_5.$ In fact if $A$ goes to
$e^{\pm \frac{\pi i}{5}},$  then  $\l$ goes to approximately 13.4721. If $A$
goes to $e^{\pm \frac{3 \pi i}{5}},$ then  $\l$ goes to approximately 4.52786.
Let $\k$ be the positive square root of $\l$. Then the roots  of  $\G_{5}
\left(D_0(F8)\right)$ are $\frac{1+ \k i}{2}, $ and $\frac{1- \k i}{2}.$ So we
have:

$$<(S^3(D_0(F8)))_d>_5 =   \left(\frac{A^{2d}}{2^d}\right)
((1+ \k i)^d + (1- \k i)^d) =
\frac{A^{2d}}{2^{d-1}}\sum_{r=0}^{[d/2]} (-1)^r \binom{d}{2r} \l^r \tag 6.5$$
The behavior of the argument (or phase) of $<(S^3(D_0(F8)))_d>_5$ mod $\pi$ as
a function of $d$ is quite simple in this case.

\subhead The covers of 0-surgery along the knot $8_1$\endsubhead

This example is more typical. The 3-twisted double of the unknot is
$8_1.$ One may of course easily compute $<(S^3(8_1))_d>_5$
exactly by recursion. For example
$<(S^3(8_1))_{17}>_5=  188 + 152 A + 136 A^2.$  There is not much pattern.
However if we embedd $k_5$ in
$\Bbb C$ by sending $A$ to $e^{\frac{\pi i}{5}}$ then the eigenvalues of
$Z_5(8_1)$ are $e_1\approx 0.676766 - 1.2548 i,$ and $e_2\approx .632251 +
0.303739i.$
$e_1$ has norm greater than one and $e_2$ has norm less than one.
Thus $<(S^3(8_1))_d>_5=(e_1)^d +(e_2)^d.$
So ${<(S^3(8_1))_d>_5}_{|A=e^{\frac{\pi i}{5}}} \approx (e_1)^d,$ for $d$
large.

\head \S 7  A connected sum formula \endhead

If $c$ is a  $q$-color, let $K(c)$ denote
$S(K)$ with the image of a meridian colored $c.$  Here and below the banding
on a meridian is taken to be that given by another nearby meridian.  As before
we have a generator $\c \in H^1(S(K)).$
Let $Z_p(K,i)= Z_p((K,i),\c).$ Define $\G_p$  similarly.  Note $Z_p(K,0) =
Z_p(K).$
Since $V_p$ of a surface with a single odd colored banded point is zero,

\proclaim{Proposition (7.1)} If $c$ is odd,  $Z_p(K,c)=0$ \endproclaim

Because $V_p$ of a 2-sphere with a single colored banded point is zero, we
have:

\proclaim{Proposition (7.2)} If $c \ne 0,$ $Z_p(U,c)=0$ \endproclaim

Let $U(i,j,k)$ denote  zero surgery to the unknot with the images of three
meridians colored $i,$ $j,$ and $k,$ where these are are good $q$-colors.  By
\cite{BHMV1,(4.4)}, $V_p$ of a 2-sphere with three points colored
$i,$ $j,$ and $k,$ is one dimensional if $(i,j,k)$ is  a small admissible
triple, and is zero otherwise.
Thus

\proclaim{Proposition (7.3)} Assume $p\ge 3.$  $Z_p(U,i,j,k)$ is the identity
on a free rank one $k_p$-module if $(i,j,k)$ is  a small admissible triple, and
is zero otherwise.
 \endproclaim

By (7.3), and  the Colored Splitting Theorem \cite{BHMV1,(1.14)}, one has:

\proclaim{Theorem 7.4} Assume $p\ge 3.$
$$Z_p(K_1\# K_2)= \bigoplus_{i \text{\ is a $q$-color}}
Z_p(K_1,i)\otimes Z_p(K_2,i).$$
 More generally:
$$Z_p(K_1\# K_2,i)= \bigoplus_{
\text{$j$,$k$ such that $(i,j,k)$ is a small admissible triple} }
 Z_p(K_1,j)\otimes Z_p(K_2,k).$$
\endproclaim

\proclaim{Proposition(7.5)}  $Z_5(D_k(U),2)$ is periodic in $k$ with period
five. In particular, $$Z_5(D_{5n}(U),2)= [0]_{\natural}$$
$$Z_5(D_{5n+1}(U),2)= [1]_{\natural}$$
$$Z_5(D_{5n+2}(U),2)= [1-A^3]_{\natural}$$
$$Z_5(D_{5n+3}(U),2)= [1-A-A^3]_{\natural}$$
$$Z_5(D_{5n+4}(U),2)= [-A^2]_{\natural}$$ \endproclaim

\demo{Proof}  One may use the same basic procedure described in \S 5 to
calculate the colored invariants of a knot. One only needs to add a straight
colored line to the slab.
$V_5$ of a torus with one framed point colored $2$ is one dimensional
\cite{BHMV1,(4.14)}. A generator is pictured in
Figure 7. To calculate $Z_5(D_k(U),2)$ one should attach the solid handlebody
of Figure 7 to the bottom of the slab of Figure 3 with a  vertical line colored
$2$ added so that the arcs labelled $2$ match up.
This new picture represents some multiple of the generator pictured in Figure
7. Just as in \S 5 we must  multiply by $\k^{-3}$ to correct for the
$p_1$-structure. Recall $\o = \eta(1+\d z).$  But this diagram  with the curve
labelled $\o$ deleted will represent  zero since $f_2$ times  $\supset \subset$
is zero in the Temperley Lieb algebra.  Thus $Z_5(E)$ is multiplication by
$\eta \delta \mu^{2k+1}$ times $a_k,$ where $a_k$ is  the multiple of Figure 7
represented by our picture with curve label $\o$ colored $1$ and the $2k+1$
twists deleted.

One may calculate $a_k $ doing a standard Kauffman bracket calculation of the
link in the diagram colored $1,$ discarding any terms in an expansion where the
segment labelled $2$ is joined to a loop which is inessential (again
since $f_2$ times  $\supset \subset$ is zero.) One is left with $a_k $ times
the generator. One has $a_{k}= A^2 a_{k-1} + (1+A) \mu^{2-2k},$ and $a_0=0$.
The result follows easily. \qed\enddemo
$$\epsffile{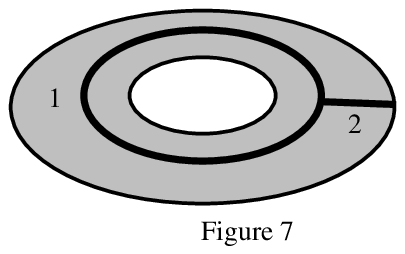}$$

The above result may also be derived from (7.7),(7.9) and (7.10). It also
follows from (8.6) below.
The above proof is more elementary and helps us understand these invariants
concretely.
By (5.10),(7.4) and (7.5),  we have for instance:

$$\align &\Cal E _5(RT \# LT )= \{1,1,1,A_3^2 , \bar A_3^2  \}\\
& \Cal E _5(F8\# F8)= \{1,1,1, A^2, \bar A^2\}\\
& \Cal E _5(RT \# RT )= \{A^4,\bar A^2,\bar A^2,A_3^2 \bar A^2,
\bar A_3^2  \bar A^2\} \tag 7.6\endalign$$

\noindent Before we made this calculation, we had found $\Cal E _5(F8\# F8)$
using the method used in calculating $\G_5(8_8).$

Let $\Cal S(c,p)$ be the set of good q-colors $i$ such that $(i,i,c)$ form a
small admissible triple.  The colored graphs in the solid torus pictured in
Figure 7 with 1 replaced by $i$ and $2$ replaced by $c$ as $i$ ranges over
$\Cal S(c,p)$ forms a basis for $V_p$ of a torus with a single point colored
$c,$ \cite{BHMV1,(4.11)}. Let $n(c,p)$ be the cardinality of $\Cal S(c,p).$ We
form two
$n(c,p)\times n(c,p)$ matrices $\eurm B(J,k,c)$ and $\eurm L(J,c)$ with rows
and columns indexed by $\Cal S(c,p)$. At this point, we begin to suppress the
dependence on $p$ again.  $\eurm
B(J,k,c)_{i,j}= \b\sum_{s=0}^{n-1}<e_s> \mu(s)^{2k+1} b(J,k,c)_{i,s,j}$
where $b(J,k,c)_{i,s,j}$ is the bracket polynomial of the colored banded link
in Figure 8a.
$$\epsffile{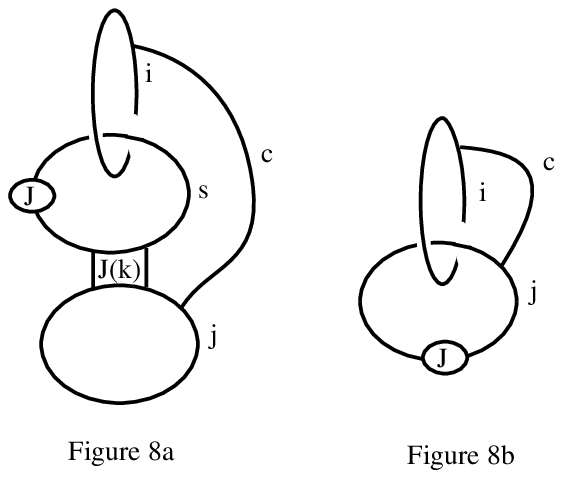}$$
 Let $\eurm L(J,c)_{i,j}$ be the bracket
polynomial of the colored banded link in
Figure 8b. Let $J(e)$ denote the knot $J$ colored $e.$

\proclaim{Theorem (7.7)} If $p\ge 3$, and $<J_{e}>_p$ is nonzero
for all $e \in \Cal S(c,p),$ then $Z_p(D_k(J),c)=(\eurm B(J,k,c) \eurm
L(J,c)^{-1})_\natural.$\endproclaim

 By the same arguments as used for (5.14), we have :

\proclaim{Corollary (7.8)} If $p\ge 3$,   and $<J_{e}>_p$ is nonzero
for all $e \in \Cal S(c,p),$  then $Z_p(D_k(J),c)$ is periodic in $k$ with
period $p.$\endproclaim

Using the tetrahedron coefficient $\left< \matrix A&B&E\\D&C&F \endmatrix
\right>,$ in the notation of \cite{MV2} ,  we have:

$$\eightpoint
\align
&\eurm L(U,c)_{i,j}=  \sum_{
\text{$r$  $\ni$ $(i,r,j)$ is a small admissible triple} } \frac{\d(r;i,j)^2
<r>}{<r,i,j>}\left< \matrix c&j&j\\r&i&i \endmatrix   \right>.\tag 7.9\\
&\text{If $(c,s,s)$ is not a small admissible triple,\ } b(U,k,c)_{i,s,j}=0
\text{,\ otherwise\ } b(U,k,c)_{i,s,j}=\tag 7.10 \\
\ &\\
 &\sum_{\matrix
\text{$r$,$r'$ $\ni$ $(i,r,s),$ and $(j,r',s)$}\\
\text{ are small admissible triples}\endmatrix }
\frac{\d(r;i,s)^2 <r> \d(r';j,s)^{2k}<r'>}{<r,i,s> <r',j,s> <c,s,s>}
\left< \matrix c&i&i\\r&s&s \endmatrix   \right>
\left< \matrix c&j&j\\r'&s&s \endmatrix   \right>.
\endalign
$$

\head \S 8  Quantum invariants of branched cyclic covers of knots  \endhead

Let $s_d(K,c)_p$ denote the sum of the $d$th powers of the roots of
$\G_p(K,c).$ Note that in \S 6, we calculated $s_d(K,0)_p$ in several cases.
The same methods, originally discussed in the remarks following (1.8),  may be
applied to calculate $s_d(K,c)_p,$ from $\G(K,c).$

\proclaim{Proposition (8.1)} If $p \ge 3$  and $c$ is even,
$s_d(K,c)_p= <K(c)_d>_p.$ \endproclaim

Let $K_d$ denote the branched cyclic cover of a knot $K$ equipped
with a \ps $\a$ with $\si(\a)= 3\si_d(K).$ This is in the homotopy class of
\pss which extend across the branched cover of $D^4$ along a pushed in Seifert
surface. Recall $<e_i>= (-1)^i \frac{A^{2i+2} -A^{-2i-2}}{A^{2} -A^{-2}}.$

\proclaim{Theorem (8.2)} For $p\ge 3,$
$$<K_d>_p
= \cases
\eta \sum_{i=0}^{(p-3)/2}<e_{2i}> s_d(K,2i)_p & \text{if $p$ is odd}\\
\eta \sum_{i=0}^{[p/4]-1}<e_{2i}> s_d(K,2i)_p & \text{if $p$ is
even.}\endcases$$
\endproclaim

 \demo{Proof}
Note that $0$-framed surgery  to  $S^3(K)_d$ along  the inverse image of the
meridian is actually $K_d.$ The trace of this surgery has zero signature, so
the $\si$ invariant of the \ps is the same for  $K_d$ and $S^3(K)_d.$

Suppose we have a surgery description of $K$ in $S^3$.
Consider the resulting surgery description of $S^3(K)_d.$   Let $\Cal D_c$ be
the result of replacing each surgery curve by $\o$ with the given framing and
then adding to the resulting picture the inverse image of the banded meridian
colored $c.$  Let $\Cal D_{\o}$ be the result of replacing each surgery curve
by $\o$ with the given framing and replacing  the inverse image of the banded
meridian  by $\o.$
We have that  $ <K(c)_d>_p= \eta <\Cal D_c>,$ and $ <K_d>_p
= \eta <\Cal D_w>.$ Thus by (8.1) if $p \ge 3$  and $c$ is even, $<\Cal D_c>=
\eta^{-1} s_d(K,c).$  Here $<X > $ just the Kauffman bracket of the linear
combination of framed links $X$, after letting $A=A_p.$
Of course instead of replacing a curve by $\o$, we could replace it by any
other combination of framed links in $S^1 \times B^2$ which represents the same
element of $V_p(S^1 \times S^1).$

 If  $p$ is even, and $c$ is odd, then by (3.1)
$<K(c)_d>_p=0.$  If $p$ is even, we are done  since $\o =\eta \O,$ and $\O =
\sum_{i=0}^{(p-4)/2} <e_i> e_i.$  For $p$ odd, we replace
$\o$ by $\o' =\eta \sum_{i=0}^{(p-3)/2} <e_{2i}> e_{2i}.$ One uses
\cite{BHMV2,6.3(iii)}, to see that $\o$ and $\o'$ represent the same element in
$V_p$ of the boundary of a solid torus.\qed \enddemo

If $K$ is a knot in $S^3$, then $K_1$ is  $S^3$. Thus we have the following
restriction on the colored invariants of $K.$

\proclaim{Corollary (8.3)} For $p\ge 3,$
$$1=\cases
 \sum_{i=0}^{(p-3)/2}<e_{2i}> \text{Trace\ }Z_p(K,2i) & \text{if $p$ is odd}\\
\sum_{i=0}^{[p/4]-1}<e_{2i}> \text{Trace\ }Z_p(K,2i) & \text{if $p$ is
even.}\endcases $$
\endproclaim

Since $\eta_3 =-1,$ we have the following corollary which may also be derived
from  (5.3).
\proclaim{Corollary (8.4)} $<K_d>_3 =-1$
\endproclaim

\proclaim{Corollary (8.5)} $<K_d>_6 =\eta_6 <S^3(K)_d>_6,$ and
$<K_d>_2 =\eta_2 <S^3(K)_d>_2.$ In particular, $<(D_k(J))_d>_2 =\eta_2 s(d,k).$
Thus $\eta_{10}^{-1}<(D_k(J))_d>_{10}= s(d,k) j_5(\eta_5^{-1} <(D_k(J))_d>_5).$
\endproclaim
\demo{Proof} The first equation is just (8.2) for $p=6.$ Applying \cite{BHMV1,
(1.5)} to $S^3$ with $\si(\a)=0$ we see $i_p(\eta_2) j_p(\eta_p) = \eta_{2p}.$
The second equation follows from the first,  $i_3(\eta_2)  = -\eta_{6},$ and
(8.4). The third equation follows from the second and (6.3).
The last equation follows from the third and \cite{BHMV1, (1.5)} again.\qed
\enddemo

We may use (8.3) to obtain a generalization of (7.5). We note that when $p$ is
five, $<e_2>^{-1}=-(A+\bar A).$

\proclaim{Corollary (8.6)}
$Z_5(D_k(J),2)= [(A+\bar A)(\b (1+\mu^{2k+1})\eurm b_k(J) -1)]_\natural.$
In particular,
$Z_5(D_k(J),2)$ is periodic in $k$ with period five. \endproclaim

For instance $Z_5(D_0(F8),2)= [A^3+A^4]_\natural.$ In general,
it is now an easy matter to calculate $<D_k(J)_d>_5$ recursively once we know
$<J>,$ and $[[J]].$ But these are just two values of the Kauffman polynomial,
for which extensive tables exist.  Moreover $<D_k(J)_d>_5$ is  periodic in $k$
with period five.  We note that $D_2(U$) is the stevedore's knot $6_1.$ We have
for instance:

$$\align
&\eta^{-1}<RT_d>_5=  (A_3 \bar A)^d+ (\bar A_3\bar A)^d+
(-1)^d(1-A+A^4)A^{2d}\\
&\eta^{-1}<F8_d>_5=  A^d+ \bar A^d+ (1-A+A^4)\\
&\eta^{-1}<(6_1)_d>_5=  1 +\bar A^d  + (1-A+A^4)(1-A^3)^d\\
&\eta^{-1}<(8_1)_{17}>_5=  1175 + 762 A + 1123 A^2\\
&\eta^{-1}<(D_0(F8))_{d}>_5=  (1-A+A^4)(A^3+A^4)^d +
\frac{A^{2d}}{2^{d-1}}\sum_{r=0}^{[d/2]} (-1)^r \binom{d}{2r} \l^r
.\endalign$$

Here $\l = 12 (A + \bar A) -8 (A^2 + \bar A^2)-1,$ as in (6.5).
Note that $<F8_d>_5$ is periodic in $d$ with period ten.
Also $<RT_d>_5$ is periodic in $d$ with period thirty. This second periodicity
is generalized below. Let $\theta_p$ denote the invariant of oriented closed
3-manifolds defined in \cite{BHMV2}.

\proclaim{Theorem (8.7)} Suppose
 $K$ is a fibered knot with a periodic monodromy of order $s$.
Then $<K_d>_2$ is periodic in $d$ with period $s.$ If $p\ge 3,$ then
$<K_d>_p,$ and $\theta_p(K_d)$ are periodic in $d$ with period $ps.$
If $r\ge 3,$ $\tau_r(K_d)$ is periodic in $r$ with period $2rs.$
 If $p=2r$ with $r$ odd and $r \ge 3,$ then
$<K_d>_p$ and $\theta_p(K_d)$ are periodic in $d$ with period $rs.$
If $r$ is odd and $r\ge 3,$ then $\tau_r(K_d)$ is periodic in $r$ with period
$rs.$
\endproclaim
\demo{Proof} The diffeomorphism type of $S^3(K)_d$ is periodic in $d$ with
period $s.$  It follows that  the first Betti number $b_1(<K_d>)=b_1(S^3(K)_d)$
is periodic with period $s.$ By \cite{DK,(5.2)},
$\si_{s+k}(K) = \si_{k}(K)+\si_{s+1}(K).$ By Lemma (5.5), $\si_{s+1}(K)\equiv 0
\pmod{8},$ so $\k^{9\si_{s+1}(K)}=u^{\frac{3}{2}\si_{s+1}(K)}$ is also a $p$th
root of unity.  Thus $\k^{3\si<K_d>}$ is periodic in $d$ with period $ps.$ In
fact if
$p$ is even, $\k^{3\si<K_d>}$ has period $\frac{ps}{2}.$
The first statement then follows  from  (8.5).

Now suppose $p\ge 3.$ By (8.2), $<K_d>_p,$ will be periodic in $d$ with period
$ps$  if  the roots $\G(K,2i)_p$ are all $ps$th-roots of unity.
The underlying manifold of $K(2i)$ is the mapping torus of the closed off
monodromy on the fiber capped off, which we denote $\Si$. $\Si$ has been given
the structure of a banded point colored $2i,$ and $K(2i)$  has been given the
extra structure of a meridian colored $2i$ with a certain banding.
Let $E$ be the associated fundamental domain of the infinite cyclic cover of
$K(2i)$
and $E_s = \cup_{0\le i \le s-1}T^i E$ as in \S 1.
$(Z_p(E))^s=Z_p(E_s)$ is then given by a scalar multiple of the identity. This
is
because
$E_s$ is diffeomorphic to $\Si\times[0,1],$ forgetting extra structure. Note
that $K(2i)_s$  has an induced \ps with $\si$ equal to $3\si_d(K).$ Whereas if
we gave $\Si \times[0,1]$ the product $p_1$-structure, then
$Z_p(\Si\times[0,1])$ would be the identity, and the associated mapping torus
has a \ps with $\si$ equal to zero.  See the proof of (6.10).  Also the banding
on the links in $E_s$ and $K(2i)_s$ differs by some number $b$ of twists. Thus
$(Z_p(E))^s$ is $\k^{9\si_{s}(K)} \mu(2i)^b$ times the identity. Note $\mu(2i)$
is a $p$th-root of unity. By lemma (5.5), $\si_s(K)\equiv 0 \pmod{8}.$ So
$\k^{9\si_s(K)}=(u^2)^{3\si_s(K)/4}$ is also a $p$th root of unity. It follows
that  $(Z_p(E))^{ps}$ is the identity and all the roots of $\G(K,2i)_p$ are
$ps$th-roots of unity.

By \cite{BHMV1,\S 2}, we have that $\theta_p(K_d)$ are periodic in $d$ with
period $ps.$ By \cite{BHMV3,(2.2)}, $\tau_r(K_d)$ is periodic in $r$ with
period $2rs.$

If $p=2r$ where $r$ is odd use \cite{BHMV1,(1.5)}
to express $<K_d>_p$ in terms of $<K_d>_{r}$ and $<K_d>_2.$ This gives the
above periodicity for $<K_d>_p.$  The periodicity of $b_1(<K_d>_p)$,
$\k^3{\si<K_d>},$ and  \cite{BHMV3,(2.2)} then yield the above periodicity of
$\theta_p(K_d)$ and $\tau_r(K_d).$
\qed\enddemo

Using Goldsmith's construction \cite{Go} of the fibration for $T(a,b)$, the
$(a,b)$ torus knot, it is easy to see that the monodromy is periodic with
period $ab.$ $K(a,b)_c$ is the the Brieskorn manifold $\Si(a,b,c)$  with a \ps
$\a$ such that $\si(\a)= 3\si_c(T(a,b).$ $\si_c(T(a,b))$ is equal to the
signature of the variety $z_0^a +z_0^b +z_0^c =1$ in $\Bbb C^3.$ There is a
well known formula due to Brieskorn  for this signature.
Our convention is that $\Si(a,b,c)$ is oriented as the boundary of the variety
$z_0^a +z_0^b +z_0^c =1$ intersected with the 6-ball. We have that
$$\align
&<\Si(a,b,c)>_p = <\Si(a,b,c+pab)>_p\\
&\t_r(\Si(a,b,c)) = \t_r(\Si(a,b,c+2rab)) \text{\ if $r$ is even} \tag 8.8\\
&<\Si(a,b,c)>_{2r} = <\Si(a,b,c+rab)>_{2r} \text{\ if $r$ is odd}\\
&\t_r(\Si(a,b,c)) = \t_r(\Si(a,b,c+rab)) \text{\ if $r$ is odd} \tag 8.9
\endalign $$

Making use of the fact $\Si(a,b,-c) =-\Si(a,b,c),$ one sees that these
equations hold for all $a,b,c$ positive or negative. In this way one obtains
four further relations,  for example:
$$<\Si(a,b,c)>_p = \overline{<\Si(a,b,pab-c)>_p}.$$
The fact that  the diffeomorphism type of $\Si(a,b,c)$ is invariant under
permutations of  $(a,b,c)$ leads to further relations among $<\Si(a,b,c)>_p.$

Freed and Gompf showed, for $c$ of the form $6 k \pm1,$ that
$\tau_r(\Si(2,3,c))$ was periodic in $c$  with period $6r$ in the case $r$ odd
and with period $3r$ in the  case $r$ even.  They used the fact that $\Si(2,3,6
k\pm 1)$ is
$(\mp 1/k)$-surgery on $\pm T(2,3)$ and   the periodicity of $\t_r$ for
$(1/n)$-surgery on a knot  due to Kirby and Melvin. In (8.9) above, we
generalize this  periodicity for $r$ odd. In (8.8), we get periodicity with
four times this period. However we have no restriction on $c$ in either case.

Let $v_p(g,c)$ be the dimension of $V_p$ of a surface of genus $g$ with a
single point colored $c.$   We may use (3.8),(8.2), and the triangle inequality
to obtain:
 \proclaim{Theorem(8.10)} If $K$ is a fibered knot with genus g, then for all
$d$

$$|\eurm i \left(\eta^{-1} (<K_d>_p)\right)| \le
\cases
 \sum_{i=0}^{(p-3)/2}\eurm i <e_{2i}> v_p(g,2i) & \text{if $p$ is odd}\\
 \sum_{i=0}^{[p/4]-1}\eurm i <e_{2i}> v_p(g,2i) & \text{if $p$ is
even.}\endcases$$

\endproclaim

\head \S 9 Partial Invariance under skein equivalence    \endhead

The invariants which we discuss are almost invariant under skein equivalence.
To see this, we first  define
some cobordism categories with more morphisms.
We define $\check C_2^{p_1}$ to be the cobordism category with objects
those of $C_2^{p_1}$ but whose morphisms are $k_p$-linear combinations of
morphisms of $C_2^{p_1}$ between the same objects.  Composition is defined as
follows. Suppose $\{C_i\} $ is a finite set of  morphisms from $\Si_1$ to
$\Si_2$
and $\{C'_j\} $ is a finite set of  morphisms from $\Si_2$ to $\Si_3,$ then
$\{C_i\cup_{\Si_2} C'_j\} $ is a finite set of  morphisms from $\Si_1$ to
$\Si_3.$
Define $$\sum_i a_i C_i \circ \sum_j a_j C'_j  = \sum_{i,j} a_i b_j\
(C_i\cup_{\Si_2} C'_j).$$

We define $Z_p$ on $\check C_2^{p_1}$ by extending linearly.    One may define
a expansion functor from
$C_{2,q}^{p_1,c}$ to  $\check C_2^{p_1},$ following the recipe in \cite{BHMV1}.
Similarly we define $\check C_{2,q}^{p_1,c}$ to be the cobordism category with
objects
those of $C_{2,q}^{p_1,c}$ but whose morphisms are $k_p$-linear combinations of
morphism of $C_{2,q}^{p_1,c}$ between the same objects.  We may define
$Z_p(M,\c)$ if $M$ is a morphism in these new categories from the $\emptyset$
to $\emptyset$ just as in \S 3.
Let $\Cal B$ be a basis for $V(\Si),$ let $\Bbb M(M,\Si,\Cal B)$ denote the
matrix for $\k^{-\si (\a(M))}Z_p(E(\Si)$ with respect to $\Cal B.$ Thus
$Z_p(M,\c)=\left(\Bbb M(M,\Si,\Cal B)\right)_\natural.$

Given a colored graph $G$ in $M$ transverse to a choice of $\Si$ dual to $\c$,
we may define the skein equivalence class of $G$  modulo $\Si$
to be the equivalence relation generated by ambient isotopies of $G$, which are
the identity in a neighborhood of $\Si$
and  the local moves  described in \cite{MV2,\S2} which takes place in the
complement of this neighborhood. If  $G$ is really a link $L$ (colored $1$),
one takes the usual Kauffman skein relation.  We clearly have:

\proclaim{Proposition (9.1)} Let $\Si$ be a fixed surface dual to $\c,$ and
$\Cal B$ is a basis for $V_p(\Si),$
  $\Bbb M(M,\Si,\Cal B)$  is an invariant of the
skein equivalence class of  $G$ modulo $\Si$. Thus, if  $\Si$ is some fixed
surface dual to $\c,$ $Z_p(M,\c)$ is well defined  on the skein class of $G$
modulo $\Si.$\endproclaim

In a similar way, we  have in the notation of \S 4:
\proclaim{Proposition (9.2)}
  $\Cal Q(\Cal T)$ is an invariant of the
skein equivalence class of  $\Cal T$. \endproclaim

\head \S 10  $V_p((S^2,m)\coprod (S^2,1 ))$ and odd Links
in $S^1  \times S^2$    \endhead

 Here we discuss another approach to odd links. We let $(S^2,m)$ demote a
2-sphere with $m$ banded points (colored one).  $K(S^2,m)$ is zero
for $m$ odd, so $V(S^2,m)$  vanishes as well. Whenever
$V_p((S^2,m_1)\coprod (S^2,m_2 ))$   is  nonzero, one may
consider $z_p(L,m_2)= Z_p\left( (S^1 \times S^2,L) \coprod  (S^1
\times S^2,L_{m_2})\right) $, where $L_{m}$ denotes $m$ $S^1$-
factors.  As a first case,  we   consider $z_p(L,1).$  It turns
out that the few invariants we get in this case are very
trivial. However our  calculation of $V_p((S^2,m_1)\coprod
(S^2,1 ))$ may be of interest.

\proclaim{Lemma (10.1)} For $p \ne 3$ or $1,$ $V_p((S^2,m
)\coprod(S^2,1 ))=0$. $V_3((S^2,m )\coprod (S^2,1 ))=k_p$ if $m$
is odd and is zero if m is even.
$\dim( V_1((S^2,m )\coprod (S^2,1)) ) = c(\frac{m+1}2)$
if $m$ is odd and is  zero if m is even.
\endproclaim

\demo {Proof}  For $p$ even and $p$ greater than two,  this
follows since  $V_p((S^2,1 ))=0$ and the tensor product axiom
(M) holds \cite  {BHMV1,(1.10)}.  According to
\cite{BHMV1,(1.9)}, there is an epimorphism  $\e_{(I\times
S^2,m,1 )}:K(I\times S^2,m, 1  ) \rightarrow  V_p((S^2,m
))\coprod(S^2,1  )))$. Here $(I\times S^2, m, 1  )$ denotes
$I\times S^2$  with $m $ framed points on $\{0\} \times S^2$
and $1 $ framed points on  $\{1\} \times S^2$. $K(I\times
S^2,m ,1  )$ is trivial if m is even, this  proves the result if
$m$ is even. From now on we assume that $m$ is odd  and $p$
is either odd or equal to two.

Now $K(I\times S^2,m ,1  )$ also has a basis consisting of the
$c(m +1 )$ diagrams $\Cal E_i$ in $I  \times B^1$ with no
crossing.   We define a  bijection $f$ from this basis to  the
corresponding basis for $K(D^3,m +1 )$ by wrapping  the
segment that  meets the isolated endpoint back around ``to
the left'' so that all endpoints  are on the same side as in
Figure 9. $$\epsffile{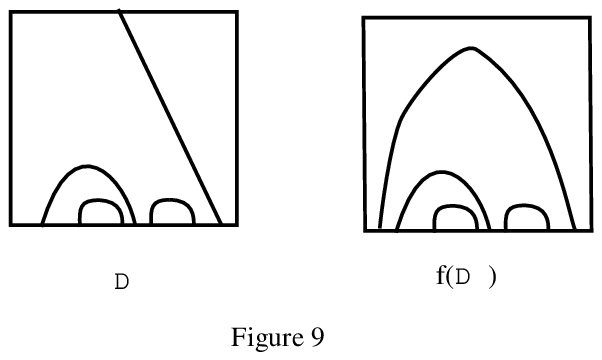}$$

As before we  define a matrix $D_p(m,1 )$ with entries $< \Cal
E_i,\Cal  E_j>_p$. This is calculated from the $<, \ >_p$
-invariant of the pair $(S^1 \times S^2 , \Cal E_i\cup-\Cal E_j)$,
where $\Cal E_i\cup-\Cal E_j$  is given by   $\Cal E_i$ glued
along $m$ framed points to the mirror image of $\Cal  E_j$
with the other points joined up in a straight fashion but
traveling  around the  $S^1$ factor in $S^1 \times S^2$. See
Figure 10 where we have drawn the  link in $S^1 \times S^2$
that we obtain when we pair the diagram of Figure  9 with
itself.

$$\epsffile{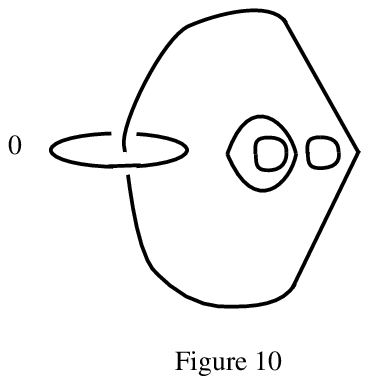}$$

We may evaluate $<\, \ >_p$ of the pair by taking the bracket
polynomial  of a linear combination of banded links evaluated
at $A =A_p$. The  specific linear combination is obtained by
replacing the component labelled  zero by the linear
combination of banded links given by $\o$ as described  in
\cite{BHMV1,5.8,\S 2}. Let  $H$ be the bracket polynomial of a
Hopf link with one component replace  by $\o$. Note in our
pictures the framing of a link is the ``blackboard  framing''.
It is clear that  $< \Cal E_i,\Cal E_j>_p$ is $\frac H{\d}$  times
$f(\Cal E_i)$ paired with $f(\Cal E_j)$ evaluated at $A=A_p$.
Thus  $D_p(m,1)= \allowmathbreak (\frac H{\d}D(m+1))_p$.
Here and  above, the bracket is normalized by saying that the
bracket of the empty  link is one.

We now observe that $H_p$ is zero for $p$ not equal to
one, or three.  It will follow that $V_p((S^2,m ))-(S^2,1 )))=0$ for $p=2$ and
for $p$ odd and greater than three. This must be true by the  ``Dirac string
trick''.
See the description of this  trick \cite{K4,(p.9)}, which seems
to be a precursor to the ``light bulb  trick'' \cite{R,(p.257)}.  If
$H_p $ were nonzero, the bracket polynomials  of the unknot
with writhe zero and with writhe two would be the same
when evaluated at $A= A_p$.  This can only happen if
$A_p^{6}=1$.   Thus $p=1$ or $p=3$.

In fact,  it is easy to show using the
axioms that for any link diagram the  bracket evaluated at
$ A_3$ is one.
Direct calculation shows $H_3=2$. Note that $\d_3 =1$, so
$D(n)_3$ is a  matrix with one in every entry . It follows
$\e_{(I\times S^2,m,1 )}(\Cal E_i)$ for each $i$ is  the same
nonzero element of $V_3((S^2,m ))\coprod (S^2,1 ))).$  Direct
calculation shows $H_1= -2 \k^3$, depending on our choice of
$(\k^3)^2=  1$. Note that $\d_1 =-2$, so $D_1(m,1)=
\allowmathbreak \k^3_1  D(m+1)_1$ which is nonsingular by
\cite{KS}.\qed\enddemo

\proclaim{Corollary (of the proof) (10.2)}If m is odd, $\{\Cal
E_i\}$ is a  basis for $V_1((S^2,m )\coprod (S^2,1 ))$.\endproclaim

\noindent {\bf Remark}    According to \cite{BHMV1,(3.9)}
If $p$ is odd and $p \ge 3$,
then $V_p((S^2,p-2 )\coprod (S^2,p-2))=k_p $. For p=3, this agrees with the
above.

Given a tangle diagram $\Cal T$ with an odd number of
stands, let $\Cal  T^+$ denote $\Cal T$ with one extra vertical
straight strand on the left. It  is not hard to prove (10.3) below.
One makes use of the fact the bracket  polynomial evaluated
at  $A_1$  is just $\d_1^{\#( D)}$, where $\#( D)$
denotes the  number of components of $D.$

\proclaim{Theorem (10.3)} Suppose L is an odd link in $S^1
\times S^2$.  Then  $z_3(L,1)$ is trivial. If $\Cal T$ is a tangle
diagram for $L$, then  $z_1(L,1)= z_1(L(T^+))$. It follows that
$z_1(L,1)$ only depends on  the absolute
values of the degrees of the individual  components of
$L.$\endproclaim

\noindent {\bf Remark} Subsequent to an earlier version of this paper,
Basinyi Chimitza \cite{C} has made some further calculations for $p$ odd and $p
\ge 3$.
He has shown that $Z_p$  on $C_{2,q}^{p_1,c}$  satisfies
a generalized tensor product axiom  in the sense of Blanchet and Masbaum
\cite{BM,Ma}, with $\hat \Si$
given by a 2-sphere with a single banded point colored $p-2.$  Also he has
shown that for $\Si$ connected $V_p((\Sigma,c)\coprod \hat \Si)$ is  isomorphic
to $V_p((\Sigma,c) \# \hat \Si).$   Moreover  $V_p(\Si \coprod (S^2,m ))=0$
if the sum of the colors of the points of
$\Si$ is odd and $m$ is odd and
$m < p-2.$  Also
$\dim(V_p((S^2,p-2 )\coprod (S^2,p )))=p-2,$ and $\dim(V_p((S^2,p )\coprod
(S^2,p )))=(p-2)^2.$

\head \S 11 Comparison with other calculations \endhead

Using our exact calculations based on (5.13), we have computed the
invariant $<\ >_p$ for zero surgery to $S^3$ along RT, F8, $5_2,$ $6_1,$
$7_2,$ for $p \le 18$.
For RT, F8, we also have computed $p=19,$ and $p=20.$
Numerical approximations of these calculations for $p=2r$ agree with the values
given in
\cite{KL,Tables:4,15,43,47,62} for $r=p/2$. Using (5.2), we have used the
calculations for $p$ odd to approximate $<\ >_{2p}.$ Again the values agree
with the values given in \cite{KL} for $r=p$. According to \cite{KL}, their
values agree with those of Freed and Gompf \cite{FG} , and those of Neil
\cite{N} wherever they overlap.

$RT$ and $F8$ are genus one fibered knots. So one can  calculate $\hat Z_p(RT)$
and  $\hat Z_p(F8)$ using the representation of $SL(2,\Bbb Z)$  specified by
Witten.  However this relies on knowing the equivalence of Witten's theory and
the theory of \cite{RT,KM,L3,BHMV1}.  I am not sure that this  has been
completely established.  Also we remark that for non-fibered knots our method
may be the only systematic way to make these computations.
 We have attempted to calculate $\hat Z_p(RT)$
and  $\hat Z_p(F8)$ using
\cite{FG} and Jeffrey \cite{Je}. In order to get the results to agree with our
earlier calculations,
we needed to take an ordered basis for the homology of an oriented fiber so
that the intersection pairing of the first basis element with the second was
minus one. This is done for the left handed trefoil in \cite{FG} without
comment. Also we noted that the left hand side of equation \cite{Je,(2.22)}
should be multiplied by three so that $\d(M,\pi)$ is integral. In view of
\cite{A3} ,\cite{KM2} and the final comment in \cite{Sc}, we modified equation
\cite{Je,(4.4)} to read $\varPsi(U) = -\varPhi(U)+3 \nu.$  Here $\nu$ is in the
notation of \cite{KM2}.   We also assumed that in the translation we should
replace $e^{2 \p i/4r}$ by $-A_{2r}$ as in  \cite{L3,Prop8}. If $K$ is either
$RT$ or $F8$, let $w_r(K)$ denote $ e^{-2 \pi i\varPsi(U)} \Cal R(U),$ after we
have replaced $e^{2 \p i/4r}$by $-A_{2r},$  where $U$ is a monodromy matrix for
the fibering.  By applying the formula \cite{Je(2.7b)} for $\Cal R$, and making
use of the Gauss sum $g(p,1)$ of \cite{BHMV1,\S 2}, we obtain

$$\align
w_{r}(RT)  &=  \left(\frac{(-1)^{r+1}A^{4-r^2}}{4r} \sum_{m=1}^{4r} A^{-m^2}
\right)\lbrack
 (-A)^{-l^2} (A^{2 l j} - A^ {-2 l j})\rbrack_{j,l}\tag 11.1\\
w_{r}(F8)  &= \left( \frac{(-1)^{r+1}A^{-r^2}}{4r} \sum_{m=1}^{4r}
A^{-m^2}\right) \lbrack
 (-A)^{j^2+2 l^2} (A^{2 l j} - A^ {-2 l j})\rbrack_{j,l}\tag 11.2\endalign$$

Here $1 \le j \le r-1$ and $1 \le l \le r-1.$
The characteristic polynomials of $w_{r}(RT)$ and $w_{r}(F8)$ agree with our
own calculations of $\G_{2r}(RT)$ and $\G_{2r}(F8)$ for $3 \le r \le 10.$
We have also checked that $w_{r}(RT)$ and $w_{r}(F8)$ are periodic in this
range. Thus $w_{r}(RT)$ is similar to $\hat Z_{2r}(RT),$ and $w_{r}(F8)$ is
similar to $\hat Z_{2r}(F8)$ in this range. Of course, this is expected for all
$r.$ We note that
$w_{r}(F8)$ is similar to $ A^{-4} w_{r}(RT) \lbrack(-A)^{4
l^2}\d_j^l\rbrack_{j,l}.$

Using
\cite{BHMV3,(2.2)} and \cite{BHMV1, (1.5)},  one has for any closed 3-manifold
$M$ with \ps $\a$ that
$$\t_5(M) =  \b_{10}^{-1}
	v^{-9-3\si(\a)} <M>_{10} \ _{| A= -v^2 \text{and } \k=v^3}.\tag 11.3$$
where $\b_{10}^{-1}= \eta_{10}^{-1}\k^{-3}=-1 - A + A^2 - A^3 - A^4 + 2A^6.$
 Of course the left hand side does not depend on $\a.$

We can compare our calculation of $(\eta_5)^{-1}<RT_d>_5$ with
Freed and Gompf's calculation  of Witten's invariant of level 3  for
$\Si(2,3,c)$ and  Neil's calculation of $\t_5$ for $\Si(2,3,c).$
Freed and Gompf  calculate these for some $c= 6k+1.$ However they observe that
$\Si(2,3,6k-1)= \overline{\Si(2,3,1-6k)}.$  This together with their
periodicity and the values they calculate is enough to determine
$\t_5(\Si(2,3,6k\pm 1))$ for all
$k.$ Neil calculated $\t_5(\Si(2,3,6k\pm 1))$ for $1\le k \le 5.$
In fact using (11.3), (8.5), and (5.4), we have:
$$\t_5((D_k(J))_d)=A^{\frac{1}{2}\si_d(D_k(U))} s(k,d)
j_5(\eta_5^{-1} < (D_k(J))_d>_5)\ _{| A= -e^{2 \p i/20}}$$
It turns out our results always agree with Freed and Gompf's. Our results agree
with Neil's for $c= 6k-1$ and are the conjugates of Neil's when $c= 6k+1.$
Presumably, this is due to a different choice of orientation conventions.  Our
orientation convention is the same as Freed and Gompf's.

\head \S 12 Afterword \endhead

Here are a few reasonable conjectures:
\medskip
{\bf Conjecture 1 }$\G_p(K)$ has coefficients which are integral polynomials in
$A_p$.
\medskip
{\bf Conjecture 2} $Z_p(F8)$ is  a periodic map.
\medskip
{\bf Conjecture 3} $Z_p(RT)$ is  a periodic map with period $3p,$
for all $p$.
\medskip
{\bf Conjecture 4} $\G_p(\text{stevedore's knot})$ has one as a root.
\medskip
{\bf Conjecture 5} The degree of $\Gamma_{2r}(\text{ tweenie knot })$ is less
than $r-1.$
\medskip
Perhaps we can tackle Conjectures 2 and 3, using (11.1) and (11.2).
\medskip
Of course there are many ways that we may begin with a familiar situation in
 link theory, and obtain a colored graph in a closed 3-manifold with a one
dimensional cohomology class. For instance given a link in a homology sphere
one may take a sublink and obtain  a manifold $M$ by performing surgery to each
component of the sublink with framing minus the sum of the linking numbers with
the other components of the sublink.    Then $H^1(M)$ is free Abelian on the
components of the sublink.
The rest of the link may then be colored. In this way we obtain many
invariants. As an example: perhaps we are interested in the symmetries of a
2-component link $(K_1,K_2).$ One may do 0-framed surgery on $K_1$ and color
$K_2$, and compare this with 0-framed surgery on $K_2$ with  $K_1$ colored.
Alternatively one could form $M$ by performing surgery along both components
and then determine if
$Z_p(M,\c)=Z_p(M,s\c)$ where $s$ is an involution on $H^1(M)$ switching duals
to the meridians.

\head Appendix   \endhead
\rm{
By a list we mean an unordered finite collection with repetitions allowed.
There is a  bijection from the set of lists of  elements of $k$ to
the set of polynomials in $f[x]$  given by
$\eufm P (\{\l_1,\l_2,\dots \l_r\})= \prod_i x-\l_i$. Given two lists we may
take the list of products of pairs of entries from each list.  Corresponding to
this product on the set of lists there is a corresponding product, which we
will also denote by $\otimes.$ The following useful formulae and their
generalizations for polynomials of higher degree are easily worked out:
 $$\eightpoint (x+a_0) \otimes (x+b_0)= x-a_0 b_0\tag A.1$$
$$\eightpoint(x+a_0) \otimes (x^2+b_1x +b_0)= x^2-a_0 b_1x +a_0^2 b_0\tag A.2$$
 $$\eightpoint \multline
  (x^2+a_1x +a_0)\otimes (x^2+b_1x +b_0) = \\
 x^4 - (a_1 b_1)x^3
	+(a_0 b_1^2 + a_1^2 b_0 - 2 a_0 b_0)x^2 - (a_0 a_1 b_0 b_1)x + a_0^2
b_0^2.\endmultline \tag A.3$$

The product $\otimes$ is closely related to the product $\boxtimes$ defined in
\cite{HZ,(p.34)}.

\enddocument
\vfill\eject \end